\documentclass[a4paper,10pt]{article}
\usepackage[utf8]{inputenc}
\usepackage{amsmath,amssymb}
\title{On very regular representations in presence of index}
\author{T. Dahn}

\begin{document}

\maketitle

\begin{abstract}
This article deals with the possibility of a very regular subgroup of a Lie group, in presence
of an index figure. Further, representations that reduce the action to a very regular boundary.
\end{abstract}

\subsection{Index}

Assume $f=-x_{1}^{2} - \ldots - x_{k}^{2} + x_{k+1}^{2} + \ldots + R$. We consider surfaces $\{ f=c \}$,
for a constant c. Example: $z^{2} \pm 1=x^{2} + y^{2}$ in $\mathbf{R}^{3}$
defines a hyperboloid. Every movement on a hyperboloid can be characterized by involutions. 
Example: For index 2 (\cite{Morse38}), consider $z_{1}^{2} + z_{2}^{2} - \zeta_{1}^{2} - \zeta_{2}^{2}$, 
for instance $((z_{1}+iz_{2}),(\zeta_{1} + i \zeta_{2} )$ in $\mathbf{C}^{2}$. 
This means that  $f(z,\zeta)=f_{1}(z) - f_{1}(\zeta)$, for instance f=0 implies that $\mid z \mid = \mid \zeta \mid$.
Further $z(f)=\zeta(f)$ implies that $\mid z \mid = \mid \zeta \mid$.
Consider $\phi(z)=(\zeta)$ as change of local coordinates. Assume $\phi^{*}d U(z)=d V(z)$, then 
we have that traces (index 2) corresponds to a change of local coordinates on a planar domain.

Consider $\overline{\{ f=df=0 \}}$, that is $\{ f=c, df=0 \}$ (\cite{Nishino68}; given $\{ S_{j} \}$
first surfaces, we have that  when $\lim S_{j}$ defines a unique first surface, the limit is regular, 
otherwise the limit defines conjugated first surfaces. 
Example: $f=P/Q$, for P,Q polynomials (of degree 2), with  constant surfaces according to $P=Q$. When $P(x)=Q(y)$, 
we assume that $x \rightarrow y(x)$ is analytic. When index is 1, y=y(x) is real, 
and when index is 2, y is complex for complex x.
Assume $y(x) \not\equiv const$ is absolute continuous and Q(y) is algebraic, then Q(y)=0 implies that x is in an algebraic set,
that is $Q \circ \phi$ algebraic when $\phi$ absolute continuous. Given Q reduced, we have that f is of
slow growth.
$y(x)$ bounded when $x \rightarrow \infty$ gives continuum for $y$. 

Assume $g \in H$ (holomorphic functions) and $(d U,d U^{\bot})$ gives an analytic continuation of g, 
then index is determined by $(dU,dU^{\bot})^{\bot}$.
When $d U$ varies over G (group as with Lie), the index gives a maximal (geometric) dimension 
for the polar relative the group, that is the complement to range of the continued continuous group. Note that $If=\int f'(t)dt=\int f(0) d u$,
that is we assume $\int d U(x,y)=\int d u=1$, further $\int \sigma_{x} \frac{\delta f}{\delta u} d x + \sigma_{y} \frac{\delta f}{\delta u} d y=If$,
with  $\sigma_{x}=du/dx$.
Note that $C(G_{2})$ can in the plane be approximated by $H(K(G_{4})) \ni d I$ (\cite{Dahn23}). Assume $\Gamma_{0}=\{ u=v \}$, then
we have that  $\Gamma_{0}(G_{2})$ corresponds to a spiral in (x,y,z),
$\Gamma_{0}(K(G_{4}))$ corresponds to $\{ u_{j}=v_{j} \}$, that is a covering of traces.
Consider $\Gamma_{1}^{j}=\{  \mid u_{j}-v_{j} \mid = R \}$, for some constant $R \neq 0$. 
We assume further that one segment $\Gamma_{1}^{i}$ satisfies a strict condition, then we can analogous 
with  Parreau (\cite{Parreau51}),  define a very regular boundary $\Gamma=(\Gamma_{0},\Gamma_{1})$. As $\Gamma_{0}^{i}$ are 
rectifiable, we can determine $\Gamma_{1}^{i} \rightarrow \Gamma_{0}^{i}$ absolute continuous (\cite{Collingwood66}).
Assume $\Sigma^{4}_{1} u_{j}$ dense in $(u_{1},\ldots,u_{4})$, 
then the traces for $g \in H$ (or $H \cap \dot{B}$) are absent. Over traces the Hessian is negatively semi-definite 
where g is convex, thus the index can in this case be determined for the restriction to absolute continuous topology.

(u,v) with  $u \rightarrow 0$ and $\mid v \mid=1$ corresponds to multivalentness for continuation $u \rightarrow v$.
Separatrices (\cite{Dulac23}) have two-valued continuation, that is v is real. 
Note that when v is real the action of the movement is not affected by v=1.

Assume sng $f=\mathcal{F}$ , we then have that  UV=VU defines a disk-neighborhood according to ${}^{t} f(u,v)=f(v,u)$. 
Density for (u,v) in a domain for f, is necessary for resolution into convex components. 
Example: $\mathcal{F}$=Hyperboloid (assume not light cone), U=$U_{1}$ translation and V=$U_{2}$ rotation, 
we then have that sng  $VU \mathcal{F}$ is discrete
and $VU \mathcal{F} \cap \mathcal{F} \neq \emptyset$ locally, that is the range of VU, $R(VU)$, is not one-sided, 
but $UV \mathcal{F} \cap \mathcal{F} = \emptyset$ locally, that is $R(UV)$ is one-sided relative $\mathcal{F}$.

When (u,v) gives a symmetric continuation of for instance a measure $d \mu$ to $C^{\infty}$ and to negative exponential type close to $\mathcal{F}$, 
we say that $d \mu$ has sharp fronts close to $\mathcal{F}$ with respect to (u,v) (cf. \cite{Garding77}). We consider measures of bounded variation close to the boundary and
usually volume preserving conjugation, that is given $d U=\rho d I$ where $\rho$ is a regular function,
the corresponding conjugation is $dV=(1/\rho) dI$.
We consider uniformities, that map a disk neighborhood of a point in (x,y,z), to the representation space (u,v) on a symmetric neighborhood
of u=v.

Assume $\mathcal{F}$ corresponds to singularities to a distribution with symbol f, that is we assume sng f
is mapped onto first surfaces to the symbol. f is assumed continuous over $\mathcal{F}$. The following main result, will be discussed in what follows.

\newtheorem{prop32}{Proposition}[section]
\begin{prop32} \label{prop32}
 Given (u,v) defines a symmetric one-sided convex neighborhood of an index figure $\mathcal{F}$, that is $\mathcal{F}(u,v) \cap \mathcal{F}(x,y,z) = \emptyset$, 
we have locally a resolution of (u,v) in convex components. Given a sharp front for f with respect to (u,v) close to $\mathcal{F}$, 
we can define a ``very regular'' subgroup. 
\end{prop32}

Assume $U_{3}$ scaling of $\mathcal{F}$ (=hyperboloid, assume not light cone), further $U_{1}=$ translation, $U_{2}=$ rotation, then $\mathcal{F}$  is generated by ($U_{1}=U_{3},U_{2})$. 
Starting from a point P on $\mathcal{F}$, we can 
consider (u,v) and $u \neq v$ as continuation, that is $U_{1}U_{2}\mathcal{F}$, is ``sharp'' relative $\mathcal{F}$, further 
$U_{2}U_{1} \mathcal{F}$ is ``diffuse'' relative $\mathcal{F}$. A symmetric one-sided disk neighborhood, can be given b 
$U_{2}U_{3} \mathcal{F}$. 
When this $\mathcal{F}$ is represented as a plane in
$(u_{1},u_{2},u_{3})$, the spiral axes $u_{1}=u_{2}$, is regarded as a normal to $\mathcal{F}$. When we are considering symmetric, contingent domains, it is sufficient to consider the 
phase, where the zero space gives $\mathcal{F}$.
When f is symbol to a ps.d.o, we consider this locally, in a neighborhood of a fixed singularity in operator space.
We assume here that index is constant in a neighborhood of the singularity.
A movement that approaches a singularity in operator space, has a corresponding movement, that we here assume
approaches $\mathcal{F}$.
For the mapping $(u,v) \rightarrow (x,y,z)$ we use a homeomorphism (\cite{Brelot51}), that is further assumed
to satisfy $sgn (\varphi(u,v)(x))_{x x}-sgn (\varphi(x)_{x x}$ has constant sign locally (Hessian) (cf. \cite{Garding77}).

The symmetry for $\mathcal{F}$ over the group G relative $U_{1}=U_{3}$, defines uniformities, in the same manner for
the symmetry relative $U_{1}=U_{2}$. Given the convergence
of exponential type, it is sufficient to consider uniformities in the phase space.
Example: given a convergence according to $\parallel \varphi(u,v) \parallel^{2} + \parallel \varphi(-v,u) \parallel^{2}$
is equivalent with  $\parallel \varphi(u,0) \parallel^{2} + \parallel \varphi(0,v) \parallel^{2}$, we do not necessarily have
a nuclear space, in this case the group corresponding to $\parallel \cdot \parallel^{2}$ is
discontinuous.  

A two-sheeted hyperboloid, can be written $L(x,x)=-x_{1}^{2} + x_{2}^{2} + x_{3}^{2}=-1$. Consider for instance
$x_{1}>0$. Then $x_{1}=cosh t^{2}$ for $0 \leq t^{2} < \infty$ is minimal for $t=0$ and $L(x/t,x/t)=-1$ iff $L(x,x)=-t^{2}$.
Alternatively, we can fix the hyperboloid and scale the axes. Given $\tilde{x}_{1}$ to Hyp(t) and $x_{1}$ to Hyp(1),
we can determine $\tilde{x}_{1}-x_{1}$, that is we can determine the traces $u_{1}(t)=u_{3}(t)$. Note
that Hyp(t) can be determined as $-(x_{1} + t)^{2} + x_{2}^{2} + x_{3}^{2}$. The light cone can be treated separately,
but is left out in this article.

In the case with  $\mathcal{F}=\{ (z,\zeta) \in \mathbf{C}^{2} \quad \mid z \mid = \mid \zeta \mid \}$ (index 2),
Assume $f$ a reduced symbol outside $\mathcal{F}$, such that $f(z-\zeta)=0$ implies that $\mid z \mid=\mid \zeta \mid$ finite. Assume
g with  $fg \sim 1$, then $g \rightarrow 0$ in $\infty$ and we can assume $z-\zeta$ dense in $(z,\zeta)$
for g. Given an absolute continuous mapping (rectifiable boundaries) from zeros to $f(z,\zeta)$ to zeros to $g(z,\zeta)$,
then g can be written as 1/f outside $\mathcal{F}$ (=$\{ g=0 \}$). Assume for instance $f \in H(K(G_{4}))$ and consider the completion to $\tilde{f} \in L^{1}$
and $\mid U \mid$ corresponding to $\int f d \mid U \mid=\int f(\mid u \mid) d \mid u \mid$, 
we have $\widetilde{\mid U \mid f}=U_{1} \tilde{f}$, that is it is sufficient to consider translation.
In the same manner, we can for a conjugated movement consider $\widetilde{\mid V \mid f}=U_{2} \tilde{f}$.
Thus can we can regard $\tilde{\mathcal{F}}$ as a spiral trace. Thus $U_{1}U_{2}=U_{2}U_{1}$ is not one-sided 
relative $\tilde{\mathcal{F}}=\{ U_{1}=U_{2} \}$, that is first surfaces to $\tilde{f} \in L^{1}$.
$(U_{1},U_{3}),(U_{2},U_{3})$ define symmetric one-sided neighborhoods of $\tilde{\mathcal{F}}$ 

\subsection{Conjugation}
Assume f regular in (x,y,z). $f(u,v,w) \simeq f(v,u,w)$, defines a two-sided limit relative w. Assume $f(u,v,w) \simeq g(u,v)h(w)$,
for instance $f \simeq (Re g)h$. In particular $\lim_{u \rightarrow 0} f=\lim_{v \rightarrow 0} f$, 
that is on R(W) we have $g(u,v) \simeq g(v,u)$, for instance W corresponds to regularization to isolated singularities. 
The symmetry for g means that g can be selected as regular locally. 

Density for translates $(u,v) \rightarrow u+v$, gives a two-sided regular limit as above.
Example: assume $y/x=const$ is mapped
on $x/y=const$ bijectively in a reflection model, or $y'(x)=const$ iff $x'(y)=const$ and $y''(x)=0$ iff 
$x''(y)=0$. In absolute continuous topology, $y''(x)=0$ iff $y'(x)=const$. 

Assume a manifold given by R(d V), is intersected by a hyper-surface N(W), 
for instance with  $d W \notin G$.
Example: $V_{i}f=f$, $i =1,\ldots,k$ implies that Wf=0. 
Thus Wf=0 intersects all leaves (invariant surfaces). Otherwise the bifurcation means that W=0
is given by a r-manifold, $r<k$. Example: given $V_{i}=I$ are disjoint, they give
limited to W=0 a (polar) r-manifold.

Starting from $\{ f=c \}$ (boundary), consider $L \cap \{ f=c \} = \{ \Omega_{i} \}$ (bifurcation).
For a normal tube, we have that  L intersects every leaf, conversely L is reached from every leaf.

\newtheorem{prop20}[prop32]{Lemma}
\begin{prop20}
Assume  U,V conjugated according to $V \simeq {}^{t} U^{-1}$. We then have for movements of bounded variation, 
such that $d {}^{t} U = \rho d I$ and $d V=(1/\rho) d I$, that
the conjugation preserves invariants, when it preserves type. When $\rho(\phi(x)) \sim 1/\rho(x)$ 
defines conjugated invariants, relative a bijective involution $\phi$, this conjugation defines contingent and symmetric
domains relative $\rho(x)=\rho(\phi(x))=1$.
\end{prop20}

Example: Assume $\rho,\vartheta \in H$, $x \rho (x)=\vartheta(1/x)$ gives $\rho=0$ for $1/x,x \neq 0$ iff $\vartheta(1/x)=0$.
Assume $<dU,dV>=1$ or $< \rho d I,(1/\rho)d I>=1$ (in the phase), given a rectifiable boundary,
we have that $\{ \rho=0 \} \rightarrow \{ (1/\rho)=0 \}$ preserves sets of measure zero (\cite{Collingwood66}). 
Note that given ${}^{t}V \sim U^{-1}$ gives conjugation,
where $U,U^{-1}$ are of the same type, given $V,{}^{t}V$ are of the same type, then the 
conjugation preserves type.

Assume $f \in (I)$ with convex support, where $d U =0$ on supp f, for $d U \in G$, considered in a convex continuous topology. This defines a closed ideal.
Given supp f included in a domain of holomorphy, there is a $d V \in G$ such that $d V(x) \neq 0$ 
on supp f and defines a separating functional, $<f,d V>$ for $f \in (I)$. Thus the measure dV does not
preserve zero lines, that is we have a maximal domain for projectivity for measures, a polar for regularity defines 
a domain for regular convergence.
 
Assume $d I \notin G_{\mathcal{H}}$, where $\mathcal{H}$ is a topological vector-space, then $Uf \rightarrow f$ does not imply that $d U \rightarrow d I$. 
Consider $(d U,d V)$ not dense.
For instance $d V=\rho d U$, with $\rho$ absolute continuous, why $d \rho=0$, gives traces. However, in the plane,
$K(G_{4})$ gives a very regular approximation property, that is $\exists (d U_{i},d U_{i+1}) \rightarrow (d I,0)$ regularly. 

Consider $\{ (u,v) \quad F(u+v) < b \} \subset \subset \Omega$ and $F(u+v)=F(u) + C_{v}$, for constants $C_{v} > 0$, 
we then have  $\{ u \quad F(u)< b \} \subset \subset \Omega$. Given $\varphi(u,v)=F(u+v) - F(u)$, with  
$\varphi$ non-linear, that is multivalued in (u,v), then sub-level surfaces are not necessarily preserved. 
Consider $f(u,{}^{t} u) \rightarrow f_{0}(u+{}^{t} u,u-{}^{t} u)$. When u+v harmonic it is real, 
v real implies that $v \sim {}^{t} v$, thus $u-{}^{t} u \sim 0$, 
that is $f_{0}$ essentially of real type. .

Consider $d U=\rho dI$, with $d U^{\bot}=(1/\rho) dI$, that is volume preserving conjugation. 
Density for $(d U,d U^{\bot})$, that is $d U+d U^{\bot}=dI$, implies a trivial polar.
Given $(d U,d U^{\bot})$ gives a projective (orthogonal), but not dense decomposition of range, the spirals are polar, 
that is we do not have a regular approximation principle for the complement of range. 
$d U + d U^{\bot}=d I + d V$ implies that if $(d U + d U^{\bot})^{\bot}=(d U + d U^{\bot})$, 
that $d V=d V^{\bot}$ that is spiral. When $d U^{\bot}=\rho d U$, for $\rho \in L^{1}$, zero-lines are preserved. 
Consider pure mappings in the plane, $(U,U^{\bot})^{\bot}=(U^{\bot},-U)$.
Assume $f^{\bot}(u,u^{\bot})=f(u^{\bot},-u)$. Then an approximation property in u for $f^{\bot}$,
has a corresponding approximation property for f in $u^{\bot}$.
 
Assume (p,q) gives zero's and poles and $L(p,q)$ a chain between them. 
According to Schwartz (\cite{Schwartz52}), when $\int_{L(p,q)} \varphi=0$ for every holomorphic differential, 
there is a  meromorphic representation $\Phi$ with  zero's p and poles q. Thus, there is given an Abel conjugation
$d U \rightarrow d V$  as above, a meromorphic representation with  an approximations property in a two-mirror model.
More precisely, according to $u=0 \rightarrow v=0 \rightarrow u=0$,
that is a continuous mapping from invariants to U to invariants to V. Consider in particular conjugation 
according to $\rho \rightarrow 1/\rho$ absolute continuous.

\newtheorem{prop33}[prop32]{Lemma}
\begin{prop33}
Assume nbhd p $\subset \{ F < b \}$, with a separation axiom: $\overline{nbhd p}$ is closed in 
$\{ F < b \}$ for b sufficiently large. When the closure is taken in $C^{\infty}$, we have point-wise  convergence, 
that is G is acting algebraically in $\overline{nbhd \quad p}$. 
If further nbhd p is disk closed, we can satisfy Petrowsky's local condition (cf. \cite{Garding77}).
\end{prop33}

Note that when Petrowsky's condition is in operator space, we have to assume nbhd p preserves negative type.
Select b minimal such that nbhd p is disk closed in $\{ F < b \}$.  Example: $ F(c(x-y))  \leq b \mid x-y \mid$,
when $\mid c \mid < 1$ 
with F(0)=0. Linear (convex) continuation corresponds to $x \rightarrow \infty$, when y fixed and 
$y \rightarrow \infty$ when x fixed. Note that not all sub-level surfaces have disk closed subsets. 

Analogous to (\cite{Garding77} Theorem 3), 
when f has singularities corresponding to $\mathcal{F}$, $\int_{bd A} f-\int_{A} d f$ has sharp fronts, 
for a very regular domain (cf. \cite{Brelot51}) A to a very regular subgroup (cf. Proposition \ref{prop32}), with very regular boundary (cf. \cite{Parreau51}) bd A. A regular path approximating sing supp,
has a dual (Fourier) regular path, approximating $\mathcal{F}$. The path is determined by (u,v).
Note that if the conjugation can be given by
L(u,v)=1, a disk closed set is invariant
for $L \rightarrow {}^{t}L(u,v)=L(v,u)$. More precisely, consider for a volume preserving conjugation,
the very regular domain $\int_{A} du dv=1$, then a symmetric, regular path on this domain that approximates the singularity, satisfies Petrowsky's local condition.

Assume $d U \rightarrow d U^{\diamondsuit}$ (harmonic conjugate) according to volume preserving conjugation over movements of bounded variation.
Assume $\Phi(d U)=d V=\sigma d U^{\diamondsuit}=\sigma \rho dU$. Given $\sigma \rho$ has $T(\sigma \rho)=\int \mid (\sigma \rho)' (t) \mid d t$ 
 (\cite{Riesz56}) constant (=1),
the transformed conjugation preserves volume. We can assume $T(\sigma \rho)$ absolute continuous.
Assume $\sigma \rho$ can be decomposed in a finite number of factors, with  $T(\sigma_{j} \rho)$ constant (=1) and $\cap \mbox{supp} \sigma_{j} \rho=\{ 0 \}$, 
then the factorized conjugation preserves volume.
Given real linear (volume preserving)
conjugation in phase to f, $u \rightarrow v$, when f algebraic in u, then f is algebraic in (u,v). 

Assume $R(d U^{\bot}) \subsetneq R(d U)^{\bot}$, where $d U^{\bot}=\rho d U^{\diamondsuit}$,
thus $R(d U^{\bot})$ gives inner points for $R(d U)^{\bot}$, further when $\rho \neq 1$ and
$d \rho \neq 0$, we have that  $R(d U^{\diamondsuit})$ gives inner points to $R(d U^{\bot})$.

Assume $u \sim v$ conjugation according to $d U = \rho d  V$, with  $\rho \rightarrow 0$ in $\infty$, 
 or $\rho \rightarrow \infty$ in $\infty$ (pseudo-orthogonal), 
we then have that $\rho \rightarrow 0$ implies that $1/ \rho \rightarrow \infty$ and $\rho + 1/\rho \rightarrow \infty$
and $d U + d V=(\rho + 1/\rho) dI$ gives invertibility.
In particular, $\rho \rightarrow 1$ implies that $1/\rho \rightarrow 1$.  

Consider f(u,v,w)=f(v,u,w), where w generates the polar.
Assume $f(u,v,w)=g(u,v)h(w)$, with  g symmetric. For instance $\mid f \mid = \mid g \mid$, where 
$\mid g \mid$ is self transposed or g(u+v)=g(v+u), that is given $g \in \dot{B}(R(W))$
we can assume f symmetric in u,v. 
 Example: a rotational surface $\mathcal{F}$ is decomposable into (traces,rotation).
When $\mid f \mid=\mid g \mid$, we have that the envelop ($\mid \cdot \mid$) has loss of dimension.
For Id : (u,v) $\rightarrow$ (x,y,z) to preserve removable sets, it is necessary
that the dimension is preserved.
 
Assume W has regularizing action, further g(u,v) = g(v,u) in point-wise  topology. A maximal domain for projectivity 
corresponds to a minimal domain for R(dW) regularizing. The dimension can be determined on R(d W), that is maximal rank on R(dW) 
does not imply maximal rank on R(d I).
An approximation property for R(dW), when $d W \rightarrow d I$ is sub-nuclear, 
implies that R(dI) has nuclear topology.
Example: $\mid T \mid^{2} \in L^{1}$ implies that $T \in L^{2}$, that is an approximation property 
for $\mid f dI \mid$ does not imply an approximation property for f d I.

\newtheorem{prop12}[prop32]{Lemma}
\begin{prop12}
As long as we consider symmetric, contingent neighborhoods of invariant sets, it is sufficient to consider the phase.
Consider $<f,d U>(\phi)$, that is f(u) can be seen as linear over a pseudo base. When dU is determined through conjugation,
it is given $f=e^{\psi}$, sufficient to determine $\phi$, such that $<\psi,d V>$ and $<\phi,d {}^{t} V>$ are linear, 
for a conjugated movement dV. 
\end{prop12}

Assuming $U I(e^{\psi}) \simeq \widehat{I}( U^{\bot} \psi)$, it is sufficient to assume $\psi,\phi \in \mathcal{S}$.
Given f(u) is regarded as algebraic and given an approximation property
that preserves this property, f can be regarded as algebraic (topologically).
Assuming $\int_{A} d f=\int_{d A} f$, thus $\int_{A} e^{\psi} d \psi=\int_{dA} e^{\psi}$.
When the equality holds in absolute continuous topology on bd A and $d \psi=0$ on A, we must have f =0 on bd A, that is $d \int_{bd A}f=f=0$.
Assume $M(x,1/x)=M(1/x,x)$,
for instance $M(1/x)=1/M(x)$, with  M algebraic. $M^{2}=1$. Example: Assume $M_{1}(x)=P(1/x)$ and $M_{2}(x)=Q(1/x)$,
that is M preserves constant value in both arguments. Thus $P(1/x)Q(x)=P(x)Q(1/x)$.
Given $M_{j}$ algebraic, $Q/P =P/Q$. In the same manner, if $M_{j}(1/x)=\overline{M}_{j}(x)$, that is $\mid P \mid^{2}=\mid Q \mid^{2}$,
that is given preservation of a constant value, it is sufficient to prove convergence for $\mid \cdot \mid$.
Given an approximation property, we can also consider a non-linear 
condition for the phase, for instance multivalued conjugation.

Assume G a parametrix to F(D) and assume $E=\widehat{G}(\phi)$, with  $\phi$ selected suitable, 
that is I is taken over $\widehat{\phi}$, when we assume FE-I has negative exponential type. 
$FE-I \sim 0$ means $FG-I=0$ (modulo $C^{\infty}$).
Assume $R(d W) \subset C^{\infty}$ and of negative exponential type on a topological vector space $\mathcal{H}$,
that is $Wf \sim 0$ (modulo regularizing action) for $f \in \mathcal{H}$.
Note for a conjugation, for instance relative $FE-I \sim 0$, that $d U \rightarrow d V$ is not necessarily 
single valued in the plane. Example:
d U=d I implies that $d V=0$ multivalued, where the order is determined on R(d W).
Note that when the strict condition holds for all derivatives to FG-I ($\in \dot{B}$), it corresponds to negative
type for FE-I.

\newtheorem{prop14}[prop32]{Lemma}
\begin{prop14}
Given $FE-I \sim 0$ on R(dW), with  conjugation $UFVE - I \sim 0$ relative $\mathcal{H}$ and 
$d I \in G_{\mathcal{H}}$, this implies that ${}^{t} V \sim U^{-1}$.
That is the contra-gradient representation to U is given by the conjugated movement on the symbol to a parametrix.
When $0 \neq FE \neq$ I, we have point-wise  convergence.
Given $C_{c}^{\infty}$ approximates $\mathcal{H}$ (density), we have $V \sim {}^{t} V$, 
that is does not change type.

\end{prop14}

Assume every movement acting on E, has a corresponding movement acting on G. We use the same notation for the 
corresponding movements.
Assume FG-I $\in C^{\infty}$ that is we have point-wise  continuation over FG-I. Further,
$({}^{t} VU - I)FE \in C^{\infty}$.
Given $d \mu=d ({}^{t} VU-I) \in \mathcal{E}^{'(0)}$, we have that  $0 \rightarrow 0$ and 
sufficient for $<E(\phi),d \mu>=0$ $\forall \phi$, is that $E(\phi) \in C^{0}$.

Example: Consider dU + dV =dI + dW. Where dW=0, $d U=d I$ implies that d V=0 (conjugation relative a projective decomposition), that is $\rho d (U-I)=d V$. 
Assume further $d W=\sigma d V$, that is the union over the 
components in $R(d U)^{\bot}$ gives the support for $\sigma$.
Assume $R(d W)^{\bot}$ has conjugation relative projective decomposition, that is R(dW) is polar to projective decomposition.
Given a very regular boundary, $\sigma_{j} d V_{j}= d V_{0}$,
such that $d V_{j} \rightarrow d V_{0}$ absolute continuous, Then supp $\sigma=\cup$ supp $\sigma_{j}$.

Given $\frac{\delta^{2} f}{\delta u \delta v}=\Delta f=0$, we consider harmonic conjugation. 
For regular approximations
we have $f(0,v)=0$ implies that $\frac{\delta f}{\delta u} \neq 0$. Example: $\Delta f=0$ with  
$\frac{\delta f}{\delta v}=const \neq 0$ in (u,v), thus $f$ does not have a
regular approximation property in u, that is we have a polar relative
regularity symmetry. 

Where $\Delta f=0$, given $u \rightarrow v$ involution, there is an orthogonal resolution.
If on a domain bounded by a closed curve, $\int \Delta f d u d v=0$, then
$\Delta f=0$ for a non-real integral curve in (u,v), 
that is $v \neq 0$ outside the boundary. 
Consider $f_{u} = f_{v} + \sigma$, where $\sigma_{v} \neq 0$. When u real, $0=\Delta f=f_{uu} - \sigma_{u}$, 
that is $\sigma$ monotonously increasing in u iff  f convex in u.

\subsection{Very regular boundary}
f(u,v) has a connected domain, if a,b can be joined by a path in the domain, thus the
order of conjugation (relative the group) denotes multiple connectedness,
that is u=0 is branch-point in the domain. Sequential continuation assumes a common point u=$\mid v \mid-1$=0,
but this does not imply that $u \rightarrow v$ is single valued.
Consider $R(d U_{1}) \cup R(d U_{2})=\Omega$ as sequential movements,
then $Cl(\Omega)$ (convex closure) contains spirals as traces. 

Consider $f_{0}(x,y) \sim f(0,y)$, as a continuation of a representation of f, close to x=0. 
The continuation defines a closed curve or not. We can now represent $f_{0}$ as absolute convergent, 
in $\mid x \mid < \delta$,
$1-\delta \leq \mid y \mid \leq 1+ \delta$. The representation $f_{0}=\Sigma B_{\nu}(x)y^{\nu}$ 
has negative $\nu$,
otherwise $f_{0}$ would be unique and regular in a disk $\mid x \mid < \delta,\mid y \mid < 1+\delta$ (\cite{Hartogs06}). 
Consider $d U \rightarrow d V$ a volume preserving conjugation, such that (U+V)f(x,y)=f(u,v) and $Uf=f$ 
iff $f(0,v)(x,y)=f(x,y)$ (assume=0, 
that is $B_{\nu}=0$ $\forall \nu$), $B_{\nu}$ determines (regularity) $u \rightarrow v$,
that is we must use a discontinuous convergence in (u,v) in the latter case.

Remark 1: note that sequential approximation on one side of u=v, can be completed to include
the trace set, that is a triangle in (u,v) can correspond to a spiral in (x,y,z). 
An approximations property (dependent of topology) would imply that a closed curve approximates a 
spiral! Thus we must in this case use
a discontinuous group. 

Remark 2: Assume $\{ u=u^{\bot} \}^{\bot}$ corresponds to a projective decomposition. According to the 
parallelogram rule
$\mid d U-d U^{\bot} \mid^{2} + \mid d U + d U^{\bot} \mid^{2} = 2(\mid d U \mid^{2} + \mid d U^{\bot} \mid^{2}$).
Assume $d U = d U^{\bot}$ and $d U + d U^{\bot}=d I$, with  $d \mid  I \mid=d I_{+} + d I_{-}$, 
that is we assume $d \mid U \mid(\phi) \simeq \mid d U(\phi) \mid$.
Given the parallelogram rule holds, Banach-spaces are Hilbert-spaces, that is not necessarily nuclear, 
that is we do not necessarily have representation in a continuous group.

Consider $u \rightarrow v$ as a continuation, when $U \rightarrow I$, a projective decomposition continuation assumes $V \rightarrow 0$ regularly. 
Given a trivial trace, the continuation (u,v) is an entire line. In particular, using the Exp-mapping, 
we can map an entire line onto a closed curve.
When Exp(u,v) traces a closed curve, that is $d f(u,v)=d \widehat{g}(u,v)=d g(Exp (u,v))$, dg is closed on zero-lines to df. 
$\int_{\Gamma} d \widehat{g}(u,v)=0$ for $\widehat{g}$ analytic, 
is invariant for change of local coordinates, for lines $\Gamma$. 

Consider $\int_{(\Gamma)} d U(f)=\int_{\Gamma} Uf$.
Assume $\Gamma$ very regular as follows:  $D \subset \subset S$, where S is an open Riemann surface. 
Define $\Gamma=\{ e^{u} \}$ with  u subharmonic and
$\lim \sup_{D \ni p \rightarrow q} u(p) \leq f(q)$, for $q \in bd D$. The upper envelop is
$\underline{H}_{f}=\sup_{\Gamma} u(p)$. Obviously $\underline{H}_{f} \leq f \leq \overline{H}_{f}$
on bd D. bd D is very regular, if we have existence of an analytic segment in $\Gamma$, 
with  isolated singularities (cf. a strict condition is given in \cite{Parreau51}).
When f is continuous and locally bounded in nbhd q, 
given bd D very regular, we have $\lim_{p \rightarrow q} H_{f}=f(q)$ (\cite{Parreau51})

Note that $u(p) \rightarrow c$, can be continued harmonically to $u(p) \rightarrow 0$. Example:
$\int_{\Gamma} d \tilde{u}=\int_{\tilde{\Gamma}} d u$,
for instance $d u \rightarrow d \tilde{u}$ absolute continuous, that is $d \tilde{u}=0$ over $\Gamma$ 
implies that $d u=0$ over $\tilde{\Gamma}$,
corresponding to discontinuous convergence. Example: given an annulus $z=u-v$, then $ \mid z \mid=0$ gives 
a radius and  $\mid z \mid=R$ an outer boundary circle, that is we have an absolute continuous mapping between
the boundaries, where the radius is transversal to the circle, a strict condition. 
The two boundaries define a very regular boundary.

Assume $A=\{ g < b \}$, with  $A \cap H_{2} \neq \emptyset$ (2-dimensional web), defines traces as inner points. 
Given g regular and $A \subset \subset \Omega$, we have that  max g on bd A. 
Thus d g=0 on bd A.
Consider $\Gamma=\Gamma_{0} \cup \Gamma_{1}$, where dg=0 on $\Gamma_{1}$ implies that x=P a point 
and d g=0 on $\Gamma_{0}$ implies that $x \in L$ a line. Assume $\psi : \Gamma_{1} \rightarrow \Gamma_{0}$ 
absolute continuous. dg is closed if dg=0 on $\Gamma$, such that $\Gamma \sim 0$, this is here substituted 
with existence of $\psi$ absolute continuous.

Assume bd=$\{ u = v \}$ and $\mathcal{F}$ is 2-dimensional in $R^{3}$, such that $\mathcal{F} \cap bd =L$ a segment
(L continuum). Consider symmetric neighborhoods of the boundary, such as $f_{uv}=f_{vu}$. 
u-v=0 can be regarded as
inner points relative the convex closure. 

\newtheorem{prop1}[prop32]{Proposition}
\begin{prop1}
Assume $\Gamma=(\Gamma_{0},\Gamma_{1})$. Consider $\Omega_{0}=nbhd \Gamma_{0}$ and $\Omega_{1}=nbhd \Gamma_{1}$
Assume $d \mu=dI(u,v)$, has bounded variation, but $\notin G_{\mathcal{H}}(\Omega_{0})$, 
where $\mathcal{H}$ is a topological vector space (for instance Hilbert).
Then, given $\Gamma_{1} \rightarrow \Gamma_{0}$ absolute continuous, $d I \in G_{\mathcal{H}}(\Omega_{1})$ can be continued to $\Gamma$.
The result is independent of scaling, when the conjugation $u \rightarrow v$ is volume preserving 
and $d I \in G_{\mathcal{H}}(nbhd \Gamma)$. 
\end{prop1}

When $d I \in (G_{2})_{\mathcal{H}}$, for $\mathcal{H}=L^{1}$, the approximation principle is not dependent on both u,v. Presence
of traces implies that $d I \notin (G_{2})_{\mathcal{H}}$, but we can still have $d I \in (G_{4})_{\mathcal{H}}$. 
Example: dU+d V=dI + d P, with  $d V=d U^{\bot}$ means that where dU is reflexive, we have traces 
$d P=d P^{\bot}$.

Starting from a continuation according to $1-\delta < \mid y \mid < 1+\delta$ and 
$\Gamma_{1} \rightarrow \Gamma_{0}$ absolute continuous. When f is regular on $\Gamma_{1}$ and can be continued 
regularly over $\Gamma_{0}$, the continuation to $\Gamma$ does not contribute to index.

Given the range has sequential decomposition $R(d U) \cup R(d V)$, we have a reducible space 
, thus a resolution of Id (completely reducible cf. \cite{Garding44}). The traces define an irreducible set in the plane, 
thus absence of resolution. 

\newtheorem{prop2}[prop32]{Lemma}
\begin{prop2}
Consider $X X' \rightarrow X''$, with  $X''$ irreducible, that is maps a 
projective (orthogonal) decomposition on a decomposition with  traces. Then we have given U+V=I on 
$Y \subset X''$,$ Y \neq X''$, that $Y=0$.
\end{prop2}

Assume dW irreducible, then dW=dUV=dI implies, when for instance d W=d U, that dV=dI. Further 
$d U^{2}=d U$ implies that d U=d I. The condition $d U^{2}$ absolute continuous does not imply that d U is 
absolute continuous, but d U=0 implies that U=I. Further, $d V^{-1}U=0$  (d UV continuation) absolute continuous
implies that U=V but not necessarily that $d U =d V$, for instance when the topology is non-nuclear. 
Example: $d V=\rho d W$ irreducible, when dW maximal and irreducible, this implies that $\rho=const$. 

Consider $UVf=f(u,v)$ and $UIf=f(u,0)$. For a domain of holomorphy, we assume existence of f regular for v=0,
that can not be continued regularly through $G^{\bot}$. 
Assume $f(u,0)$ regular, when $f(u,v)$ has a root that f(u,0) does not have, we have a continuation 
principle with an approximation property, but this does not imply absence of non-regular approximations. 
We assume for an regular approximation property, that $d v/d u \rightarrow 0$ regularly.

Note that f(u,v)=g(u)h(v)=const, if $g(u)=const/h(v)$. Example: $g(u) \sim P(1/u)$ according to preservation
of constant value, with  v=1/u. 
Consider f(u,v)=0 implies that u=v=0, that is f has zero's in corners (of bi-disk). Given $\{ f < b \} \subset \subset \Omega$ and $b>0$,
f has zero's in compacts (reduced). When f(u,v)=0 implies that u=0, $\mid v \mid \leq 1$, f has zero's on edges. 
When f(u,v)=0 implies that $0 < \mid u,v \mid < 1$, f has zero's as inner points.
For traces, we note that (u,$v_{1}$), where $v_{1}$ is rotational,
defines a cylindrical spiral. $(u,v_{2})$ with  $v_{2}$ ähnlich, defines a conic spiral (cf. \cite{Lie91},
infinitesimal Spiraltransformation). 

Consider $f(\zeta)=F(\gamma)(\zeta) \rightarrow \gamma$.
Assume $\gamma$ is generated by $\gamma=\gamma(u,v)$. $\gamma$ in the 
spiral case, is dependent of both u and v.  
In the case when $u/v \rightarrow 0$ as $u \rightarrow 0$ implies that $\gamma(u,v)$ is regular, 
we can assume a  regular approximation principle for $\gamma$.
Assume Hyp = sng f. Given (u,v) is mapped on Hyp, the preimage is regular (in the sense of \cite{Brelot51}). 
Assume $\Phi^{-1}U \Phi=V$ defines conjugation, such that in weak sense $\Phi^{-1} {}^{t} \Phi =I$ and 
$d I \in G$, that is given the conjugation is sub-nuclear, a nuclear mapping assumes weakly, presence of traces. 
In the case when $d I \notin G$, the group must be extended
(or the topology), we assume $\Phi \rightarrow \Phi^{-1}$ according to $\rho \rightarrow 1/\rho$.

Given that regular approximations are defined by a domain of holomorphy, non-trivial continuations,
if they exist, are non-regular. When $Wf \rightarrow f$ in $C^{\infty}$ topology, this defines algebraic invariants. 
Example: Given a continuation that is multivalued, with  holomorphic leaves, for instance 
$d U \rightarrow d V_{j}$, $j=1,2$ and $d V_{1} \rightarrow d V_{2}$ continuous relative the same conjugation (transponation). 
Then $f(v_{1})=const,f(v_{2})=const$,  corresponding to a non-linear condition, define conjugated
first surfaces.

Presence of invariants implies absence of an approximations property (relative the movement),
thus presence of a non-regular approximations is necessary for determination of G. 
Assume $\Gamma$ denotes invariants, where $d W^{-1}U=dI$ on $\Gamma$, for instance presence of 
non-regular approximations contradicting hypoellipticity. That is convergence outside invariants 
does not imply hypoellipticity. A very regular approximation property, gives regular 
approximations on a one-sided domain (sharp fronts). When we consider Hausdorff uniformities (\cite{Bourbaki89}), singularities 
are represented by traces. Note that T(f) absolute continuous, implies that $IT(f)=T(f)$, that is UT(f)=T(f) 
iff d U=d I. Relative a subgroup $G_{2}$, we have that $0 \neq f(u,v) \rightarrow$ f(P) in H(u,v), where P=(x,y,z) a point, 
implies that f(P) is zero or regular in (u,v).

\newtheorem{prop21}[prop32]{Definition}
\begin{prop21}
 Retraction according to $d F/dx / d F/d y \simeq dy/d x \neq const$, 
is well-defined on a neighborhood of Hess=0, (degenerate points), where $dF/dy \neq 0$ (or $d F / dx \neq 0$). 
\end{prop21}

For instance, $\frac{\delta^{2} F}{\delta y \delta x}=0$, with continuation to $C^{\infty}$
on a symmetric neighborhood, where $\frac{\delta F}{\delta x} \neq \frac{\delta F}{\delta y}$ and $\frac{\delta F}{\delta x} \neq 0$ or $\frac{\delta F}{\delta y} \neq 0$.
Symmetry as in $\Sigma \frac{\delta^{2} F}{\delta y_{i} \delta x_{i}}=0$ in the
Hessianen, means that F is harmonic.
Note that $(U-I)F + IF - F=UF-F$, and given F absolute continuous, we have $IF=F$,
that is $d U=d I$ implies that U=I. 
Note that it is sufficient to consider $C^{0}$, for a continuation to $\mathcal{E}^{'(0)}$, 
that is in Hausdorff uniform topology, given retraction well-defined, we have arbitrarily close to paths, regular in the sense of \cite{Brelot51}, there are paths regular in 
the sense of \cite{Parreau51}. .

Assume $UT(\varphi)=T({}^{t} U \varphi)$, where ${}^{t} U \varphi$ is point-wise convergent. 
I is very regular in $\mathcal{E}$, but this
does not imply that I is very regular in $\mathcal{E}'$,  $I - C^{\infty}=\tilde{I} : \mathcal{D}' \rightarrow \mathcal{D}^{'F}$,
that is we have a non-trivial kernel in $\mathcal{D}'$. $\tilde{I}$ is surjective on $\mathcal{E}$, 
but ${}^{t} \tilde{I}$ is not surjective 
on $\mathcal{E}'$, possibly 1-1.

\newtheorem{prop35}[prop32]{Lemma}
\begin{prop35}
Assume (u,v) defines a rectifiable boundary $\Gamma$, with $\Gamma_{0}=\{ u=v \}$, 
$\Gamma_{1}=\{ \mid u-v \mid =R \}$, for a scalar R. Assume F regular on inner points,
given $\Gamma_{0}$ regular (as in \cite{Brelot51}) and $\Gamma_{1}$ regular (as in \cite{Parreau51}), 
then we have that  $\Omega \backslash \Gamma$ is very regular (as in \cite{Brelot51}). 
Given $\Gamma_{1} \rightarrow \Gamma_{0}$ 
absolute continuous, we have that $(\Gamma_{0},\Gamma_{1})$ is very regular (as in \cite{Parreau51}). 
Thus, a very regular distribution F defines a very regular support  
and we have monotropy with respect to paths:
that is sng of measure zero is mapped on sng of measure zero. 
\end{prop35}

Consider for instance $\phi(u,v)=u-v$ of bounded variation, this implies a determined tangent. Note
that Hausdorff uniformity, means that the domain is regular outside $\Gamma_{0}$, corresponding to 
a very regular distribution. A strict condition on $\Gamma_{1}$,
where $\Gamma_{1}$ is acting as transversal to $\Gamma_{0}$, motivates density for u-v in (u,v), 
that is a very regular domain is symmetric with respect to the trace set.

Example: consider $\{ F < b \}$, where b is the distance between points on a web. Given 
$\tilde{b}$ is the corresponding sequential distance, we must have
$b < \tilde{b}$. F assumes the minimal value on traces. 

As in \cite{Brelot51}, assume $P \in bd A$ for A an open regular set $A \subset \subset \mathbf{R}^{n}$.
Assume $\infty \notin U=nbhd P$ and $h(U) \supset V^{n-1}$ ($V^{n-1}$ a (n-1)-web), for a homeomorphy h,
such that $h(U \cap A)=\{ \varphi > 0 \}$
and $h(U' \cap A)=\{ \varphi < 0 \}$, where $U' $ is a different neighborhood. Given 
h : (u,v,w) $\rightarrow$ (x,y,z), with w fixed (rotation), we have 
$\int f(u,v,\cdot) d u dv=\int f(u,v,\cdot)(x,y,z) \rho(x,y,z) d x d y d z$.
Note that $V^{n-1}$ =web assumes Hausdorff uniform topology in a plane (u,v). We assume that when the web separates the image space, 
the condition $\{ F(u,v) < b \}$ relatively compact,
means that $F(u,v) \rightarrow 0$ implies that $u \rightarrow 0$, that is F is downward bounded in u 
outside a compact in (u,v).
Assume L: u=v a zero line to $d \mu \in \mathcal{E}^{' (0)}$, then Exp (in L) corresponds to a closed 
curve. 
When $d \mu=0$ on any line parallel to L and on one side of L,
the corresponding closed curve $\sim 0$.

\subsection{Rectifiable boundary}

Consider analytic continuation, defined by sequential continuation, that preserve dimension. Assume $d V'' =\rho d V$, where $\rho$ is linear in u,v, linear convexity implies that $v/u \rightarrow 0$. 
Assume (u,v) volume preserving and $\rho$ single valued in (u,v). Consider $L_{ac} \rightarrow L_{c}$ continuous, 
sufficient for this is $\rho$ absolute continuous, necessary is $\rho$ of bounded variation.
Example: assume $R(V) \cap R(V')$ is not a subset of $R(d V) \cup R(d V'')$, for instance a spiral between
two circles.

\newtheorem{prop7}[prop32]{Lemma}
\begin{prop7}
Assume approximations are given by $W \subset G$, such that we have existence of V with  
$W \subset V \subset G$ and such that W is closed (regularizing) in V. This can be seen as a
necessary separation condition. When W is maximal in V, 
presence of an interpolation principle implies presence of non-regular approximations. 
\end{prop7}

The implication $dW=0 \Rightarrow d V=0$, defines closedness for W in V, relative a majorization condition.
Example: $dV=\rho dI, dW=\vartheta d I$ and $d V=\sigma d W$, with  $\rho / \vartheta \rightarrow 0$ in $\infty$, 
for instance dW has negative type relative dV, that is $\sigma$ has negative type. 

Assume $W \subset V \subset G$, then W is relatively closed in V, if bd W is irreducible and removable in V, 
this does not imply that bd W is irreducible in G,
for instance $\exists d W' \in G$ regular on $d W=d I$, that is $dW' \neq d I$. 

Consider $\mathcal{E}' \rightarrow \mathcal{E}^{'(0)}$ as regularization. 
Starting with a domain of holomorphy $\Omega$, we consider $H(\Omega) \subset C(\Omega)$,
that is a neighborhood not in H. Density can be motivated by $H(G_{8}) \subset C(G_{2}$) in the plane. 
Given a strict condition, H is dense in $L^{1}$ in $\mathbf{R}^{n}$. 
The Radon-Nikodym representation (cf. \cite{Riesz56}) of  $dU=\rho d I$, can be reduced to $\rho \in L^{1}$ and non-regular 
approximations can then be given in $\mathcal{D}_{L^{1}}'$. 
Note that inclusion between weighted Sobolev spaces, assumes a strict condition.
When H is dense in $L^{1}$ we get a continuation to $H(G_{8})(\mathbf{R}^{n})$, dense in $L^{1}(G_{2}) \simeq L^{1}(\mathbf{R}^{n})$.

Example: consider $1/(\rho + 1/\rho) \rightarrow 0$, when $1/\rho \rightarrow 0$ in $\infty$.
Further, $\rho(u,v) \rightarrow 1$ iff $1/\rho(u,v) \rightarrow 1$, given $\rho,1/\rho$ regular. 
Given (u,v) volume preserving, $\rho$ can be assumed symmetrically regular, when $\rho \rightarrow {}^{t} \rho$
preserves regularity. For a very regular domain A, $\int_{A} d W(f)=\int_{bd A} Wf$, 
that is $W \mid_{A} f=\int_{bd A}Wf$, that is the action is given by restriction to the boundary. 

Consider $fg = 1 + P_{j}$, that is on the kernel to g, we have $-P_{j}=1$ (constant surfaces), 
note that a `` semi-algebraic'' polar does not exclude presence of a spiral. Continuations of negative type 
over g=0 are regular continuations and presence of traces implies
non-negative type. Starting from $d U = \rho d V$, with  $\rho \in C^{\infty}$,  given 
$\infty \in \{ \rho = 1 \}$, we have that  R(G) has a polar.

\newtheorem{prop8}[prop32]{Lemma}
\begin{prop8}
 Consider fg=1 + $P_{j}$, with  $\{ 1+ P_{j} < \lambda \}=V_{j}(\lambda)$ semi-algebraic.
Note that a linear phase in the plane, corresponds to $P_{j}$ being algebraic.
The conjugation Uf Vg=fg, can be considered as fg(u,v)=f(u)g(v), for instance ${}^{t} V=U^{-1}$ 
where the type is preserved.
\end{prop8}

Multivalentness with  holomorphic leaves, corresponds to presence of spiral approximations (\cite{Oka60}). 
For $P_{j}$ monotonous, we consider $V_{j}(\lambda)$ as being single valued. Sufficient for this is that $P_{j}$ ($P^{N}_{j}$) 
is downward bounded. Given $P_{j} \rightarrow \infty$ 
in $\infty$, the sets $\{ P_{j} < \lambda \}$ are relatively bounded sets. Assume $f(u,v)=f(u)g(1/u)$, where 
u corresponds to zero's and v to poles,
then given $u \rightarrow v$ an Abel chain (cf. \cite{Schwartz52}), f is convergent as a meromorphic function (discontinuous function in H).
Example: $A = \{ f(u)g(1/u) \leq \lambda \}$ semi-algebraic, implies that $\{ g(1/u) \leq \lambda \}$ 
semi-algebraic. 

Assume f a symbol f=P/Q, for polynomials P,Q, on sets where fg=1.
Given a rectifiable boundary, we have $f \rightarrow g$ is absolute continuous at the boundary. 
Further, when the zero sets to f are trivial (measure zero),
we have that presence of fundamental solution implies a zero-set to g. 

Consider U acting in Exp-topology and that it gives an one-sided continuation.
Given U surjective on the range, we have ${}^{t} U$ 1-1.
Given absence of surjectivity, we have presence of polar for the range, the order for the kernel to
${}^{t}U$ is dependent of index.

Example: assume $R(dI) \rightarrow_{u} R(dU) \rightarrow_{v} R(dVU)$, where v is reversible, that is $(d U,d V)$ 
is a continuation with  $d V \in G$. Example: Assume $u \rightarrow P$, with  P singular, but
$(u,v) \not\rightarrow P$, then (u,v) is closed, in this context one-sided, 
that is (u,v) outer to a neighborhood of P. 

Consider a very regular boundary, $\underline{H}=\sup_{\Gamma} u \leq f(q)$ with  $q \in bd \Omega$. Note that u 
subharmonic, implies that u does not have relative maxima (\cite{Oka60}), further when f harmonic, 
f-u is harmonic and $\neq 0$ outside $\infty$. We have a 
maximum principle on a very regular domain, when $H(p)=f(q)$ implies $p \in bd \Omega$. 
A single valued (linear) planar continuation, gives a normal extension.

Consider reduction to a cylindrical domain, for instance $\beta_{1} d U_{1} + \beta_{2} d U_{2} + \alpha d A=0$,
with  $\beta_{1}: \beta_{2}$ non-constant,
that is absence of traces (spiral). When $\alpha \neq 0$ and d A=0, then according to (\cite{Lie91}, Chapter 19,20), two 
movements are defined simultaneously.
Assume the range a surface of revolution. Given the equation defines action in the phase, 
we consider $UVW=I$, where W defines rotation. Further, when f(u,v,w)=f(v,u,w),
then UV=I reduces action to the web. 
Consider planar sections $(x,y,z) \rightarrow (x,y) \rightarrow (u,v)$ continuous, then irreducible traces are mapped onto
irreducible traces.
When we consider (u,v) $\rightarrow$ (x,y,z),
we have a multivalued mapping in the plane, that is non-linear. Given $\{ F=const \}$ represents traces, 
we have clustersets.

\newtheorem{prop5}[prop32]{Proposition}
\begin{prop5}
Singularities can be characterized by a subgroup, that gives a regularizing continuation 
(or that preserves sng). 
Conversely, given $\mathcal{F}$ is a figure for sng, then (u,v) is
regularizing, if ${}^{t}V U \mathcal{F} \cap \mathcal{F} = \emptyset$ (given ${}^{t} V=U^{-1}$ sng is preserved) 
\end{prop5}

Example: $(U+V)f -f \sim 0$ (regular), with  $U \rightarrow V$ according to Abel, giving representations of projective type, 
corresponds to a very regular approximation. Assume E symbol to a parametrix and $UFVE \sim I$, that is 
${}^{t} V U \sim I$, also when $d I \notin G$. This conjugation is volume preserving
when $d I \in G$. For approximative solutions E, it is sufficient to assume exponential type 0 in (x,y,z) (\cite{Dahn13}).
Consider analogous to \cite{Hartogs06}, $f_{0}(x,y)=f(0,y)$ as continuation close to x=0.
Consider $T_{0}(x)$ with  $\mid f_{0}'(x,y) \mid=T_{0}'(x,y)$, where T is absolute continuous. 
The continuation of $T_{0}$ corresponds to
two-valued continuation of $f_{0}$.
Example: consider $f^{*}(x,y)=\overline{f}(y,x)$ as a closed extension, the continuation 
$f_{0}(x,y,1/x,1/y)$ is dependent of $1/x,1/y$, or it would be a regular extension.
 
\newtheorem{prop6}[prop32]{Lemma}
\begin{prop6}
Assume R reflection.
Assume $fR \leq Rf$ and $f(z)f(Rz) \leq \mid f \mid^{2}(z)$, for instance
$\mid f \mid^{2}=\sigma f(z)f(Rz)$, with $1/\sigma \rightarrow 0$ in $\infty$.
In particular, strict convergence for $T(f) \rightarrow 0$, implies convergence on a very regular boundary
for f.
\end{prop6}

That is, we assume $T(f)=\sigma T(f)T(fR)$ and a very regular boundary can be given by $\{ df=0 \} \cup \{d fR=0 \}$.
Example: one-sidedness for a domain A convex in the plane, that is $z \in A$ implies that $ Rz \notin A$ and 
If(R) $\equiv 0$, that is R ``regularizing''.

Assume C is given by a first surface $f(u)=f(0)$, in $L^{1}_{ac}$, this means that u=0, that is If=f. 
Given $d U = \rho d I \neq dI$, we have that 
$\frac{d \rho}{d t} =0$ on a discrete set. 
On a principally defined set of singularities, $\rho$ is unique. Assume $d I \rightarrow d U \rightarrow d U^{\bot}$
with $d U^{\bot}/d U=\sigma$, then given $\rho \sigma$ regular and bounded, we have a rectifiable boundary.
When $(d U,d U^{\bot})$ is considered as a continuation, where $d U^{\bot}$ is regular and invertible, 
we have a transmission property.

Consider $XX' \rightarrow X''$, with  $X''$ an irreducible boundary. Consider a boundary operator d N(u,v),
very regular over traces, that is dN - d I is regular. When we have that  d N + (d I - d N) approximates d I, 
given density, we can write $d N^{\bot}=d I - d N$. Assume B :$\Omega \rightarrow bd \Omega$ continuous.
When $\{ u+v \}$ is dense and irreducible in $bd \Omega$, we have that
$bd \Omega$ is irreducible. Given further $V \subsetneq bd \Omega$ is irreducible with  $BV \subset V$,
we have that V=0.

\emph{$d U=\rho dI$ of bounded variation, given $\{ (x,y,z) \quad \rho=1 \}=\{ (x,y,z) \quad 1/\rho=1 \}$, for restriction to lines, 
we have preservation of type for $d U \rightarrow d U^{-1}$. Note that $\{ \rho \leq 1 \}$ semi-algebraic 
implies that $\{ \rho \geq 1 \}$ is
semi-algebraic, but does not imply that $1/\rho$ is algebraic. Given $\rho \rightarrow 1/\rho$ absolute continuous, 
however sets of measure zero are preserved.}

Assume $\Omega_{j}$ are convex sets j=1,2, separated by a separating functional. Sufficient for this is that
$\Omega_{1}$ has inner points $x \in L_{1}$, where $L_{1}$ is continuum line segment, for instance when 
$\Omega_{1}=\{ F < \lambda \}$ and $\Omega_{1}$ is an unbounded set.
Thus, sufficient for presence of a separating functional to convex sets $\Omega_{j}$, is a strict condition. 
Example:F=1/f, where f is has no invariant sets relative G.
Example: $f(\phi)$ absolute continuous, has a radius of convergence, given by a separating functional. 

Assume $d U=\rho d I$, with  $\rho$ holomorphic on $\Omega$. Assume $d U^{-1}=\vartheta d I$, with  $\vartheta \in L^{1}$. 
Consider 
$\overline{\{ \rho=d \rho=0 \}} \simeq \overline{\{ \vartheta=d \vartheta=0 \}}$ 
as topological closure. Then, $\rho \rightarrow \vartheta$ can preserve type
in $L^{1}$ (topological closure), but simultaneously does not preserve type in H.
Thus it is necessary with a polar relative $\Omega$ (domain of holomorphy).

$\Omega$ is very regular (\cite{Brelot51}), if $\Omega,bd \Omega$ are regular, further we assume h a 
homeomorphism that separates inner points from the boundary.
Starting from
$C(G_{2}) \supset H(G_{8})$, for $G_{2}$ the trace set generates a geodetic on the boundary to a cylinder. 
For higher index, we can discuss traces from $G_{8}$, that is we assume 
multivalentness for traces.
Among the traces, there is a trace h(geodetic(u,v)) on an index figure $V^{2}$.  We consider traces as irreducibles, 
further $F'(x,y,z) \in C^{\infty}((V^{2})^{c})$ corresponds to $F(u,v) \in C^{\infty} (u \neq v)$, 
where $V^{2}$ is an index figure. When $\mu(u,v)(x,y,z) \leq C \mu(x,y,z)$, for 
some suitable semi-norm $\mu$, negative type is preserved. A very regular approximation principle,
assumes that regularity and negative type is preserved, when $U \rightarrow I$, $V \rightarrow 0$.

Example: given $\{ g=0 \}=V^{2}$,
with subset $\{ x_{1}=\varphi(x_{2},x_{3}) \}$, where $\varphi$ is holomorphic, that is $g(I(x_{1})-\varphi)=0$. 
Every function holomorphic on $V^{2}$ is constant (\cite{Schwartz52}).  Given g=0 defines a revolution surface, 
we always have $g(x_{1}-\varphi(x_{2},x_{3}))=0$, where $\varphi$ is symmetric in $(x_{2},x_{3})$. 
Example: consider $F(x,x^{*})$ holomorphic on a domain of holomorphy, where 
$x \rightarrow x^{*}$ is conjugated according to zero's $\rightarrow$ poles
through involution, then we have $F(x,x^{*}) \sim P(x)/Q(x)$, for holomorphic functions P,Q.

Consider convexity relative two hyperplanes. Example:
$\Omega_{-}(L_{1}) \cap \Omega_{-}(L_{2})$, where $\Omega_{-}$ is one-sided relative invariant
sets $L_{j}$, the intersection defines a very regular domain. Example: $d U=d I$ iff $d U^{\bot}=0$,
that is invariants can be given as analytic sets (not uniquely). Example: $d I \notin (G_{2})_{H}$,
that is invariants are not analytic in $H(G_{2})$,
but analytic in $H(K(G_{4}))$, in the plane. 

\newtheorem{prop15}[prop32]{Lemma}
\begin{prop15}
Assume d U + d V=d I + d W, that is dW=0 implies that (dU,dV) gives a projective decomposition, over the polar R(dW) is not assumed analytic.
Invariants dU=dI for dW=0, correspond to a reduction to the boundary (here invariants). 
Assume $L_{1}$ invariants to d U and $L_{2}$ invariants to $d U^{\bot}$, then necessary for a regular 
transmission property, 
is that $d U^{\bot}$ is invertible (or d U),
for instance regular approximations through $L_{1}$ (or $L_{2}$).
\end{prop15}

A regular transmission property in this context, refers to a property that R(d U) can be regularly 
approximated by $R(d U^{\bot})$.
For saddle points (col), we have that the characteristics of bounded variation (determined tangent) 
are all separatrices (\cite{Dulac23}). 
A very regular region, can in particular be defined by saddle points. Example: assume $\Omega$ a very regular 
region, with  $\Omega_{\pm}$ on one side of a regular surface (trace). 
Consider $\int_{\Omega} df = \int_{bd \Omega}f$, where bd $\Omega$ is defined by pseudo-characteristics,
$L_{1} \cup L_{2}$ with  $L_{1} \cap L_{2} \neq \emptyset$. That is, we assume the domain decomposed in 
two leaves. 
Example: $(L_{1}^{-}) \cap (L_{2}^{+})$ gives a
quadrant of inner points to pseudo characteristics. For a domain symmetric (and contingent) with respect to 
pseudo characteristics, we can consider $(L_{1}^{-}) \rightarrow (L_{1}^{+})$ as involution.

\newtheorem{prop16}[prop32]{Lemma}
\begin{prop16}
Convexity with respect to two planes, defines a conjugation, $\Phi^{-1}_{2} U \Phi_{1} \sim V$.
Thus, $\Phi_{1} \rightarrow I$ iff $\Phi^{-1}_{2} \rightarrow I$, implies in $L_{ac}^{1}$ that $UI \sim IU$, 
that is when $d I \in G$, we have a (topologically) algebraic movement.
. 
\end{prop16}

Assume $\Omega$ one-sided relative a very regular boundary, Assume $q(x)=0$ implies that $L(tx)=t<x',x>=0$,
for t real scalar, 
that defines a line segment. A strict condition implies existence of y close to x, with  $L(y)=0$, 
$q(y) \neq 0$, that is a very regular boundary identifies an associated separating functional. 
The consequential domain is one-sided in standard way, that is convex in the plane.
Assume $X^{\bot}$ is convex with inner points relative a line, in the same manner for $X^{\bot \bot}$,
that is reflexivity gives presence of a separating functional in this case. 
Note that bounded variation implies a determined tangent, this implies a determined normal,
that defines a single valued (in the plane) conjugation,
but not unique. Example: consider $rad J'$ as  $<x_{1}',x> \ldots <x_{k}',x>=0$ over $\Omega$.

Assume $T'(x)=\mid f'(x) \mid$.  $T(f)=\int f d \mid I \mid=\int \mid f'(x) \mid d x$.
Assume  $T(f)$ is convex, that is $T''(f) \geq 0$. 
Assume $g=\overline{f}$, that is $T'' =f'' g' + f' g''  \geq 0$, since $f' g' \geq 0$, the condition is
$(\log f')' + \log (g' )' \geq 0$ and a sufficient condition
for positivity is $f''/f + g''/g > 0$.
Assume $\log f' =\phi$, then the condition is $Re \phi' \geq o$. 

We can consider a subgroup for invariance. When $UV=I$ for $d I \in G$, we have that $U=V^{-1}$ or U=V=I. 
Example: for a 
continuation relative a projective decomposition, consider
$D(G_{\mathcal{E}^{'(0)}}) = \{ dU \quad dU^{\bot} \in G_{\mathcal{E}^{'(0)}} \}$, with  $U^{\bot}$ invertible, 
since $G_{\mathcal{E}^{'(0)}} \ni dU^{\bot}=dI-dU$ is invertible, 
we have $G_{\mathcal{E}^{'(0)}} \subset D(G_{\mathcal{E}^{'(0)}})$ in this case.
 Example: $\mathcal{F}$ hyperboloid =R(V),$R(U) \cap R(V)$ 1-dimensional, defines a subgroup (rotation),
when $U=U_{1}$ is translation, we have $UV \neq I$.  Example: $R(U_{1}) \cap R(V)$ 
0-dimensional and $R(U_{2}) \cap R(V)$ 1-dimensional, that is given R(V) defines a trace, 
the dimension is dependent of the character of the movements in the group.

\subsection{Singularities}

Note that f of bounded variation implies that $f=f_{+} - f_{-}$. Given an approximation property, 
we can discuss $f_{vu} \rightarrow f_{x_{i}y_{i}}$, that is we can discuss change of sign under $U \rightarrow I$, thus the conjugation can preserve convexity for f or not. 
Example: f convex (absolute continuous)
on the domain $\Omega$ and as long as (u,v) preserves $\Omega$, as $U \rightarrow I$, convexity is preserved. 
Given $f_{uv} \rightarrow f_{vu}$ changes sign, we have saddle points. 
 
Assume $d U \rightarrow d V$ preserves bounded variation, for instance $\mid d V \mid \leq C \mid d U \mid$.
$d V/d U \rightarrow 0$ implies that $(d U,d V)$ is strictly continuous, implies inclusion 
between the respective ideals.
But given d U=d V, a non-trivial trace, where (u,v) has bounded variation, $(d U,d V) \rightarrow (x,y,z)$ 
is not necessarily continuous.

\emph{When $d U \in \mathcal{E}^{'(0)}$, we have 
$d (I-U)^{-1} \sim \Sigma d U^{j}$. Given translation is dense in the domain, iteration can be considered
as convolution. Given a projective decomposition $U+V=I$, we have that $d V^{-1} \in \mathcal{E}^{'(0)}$}

Example: $G(e^{U(x,y)})=\widehat{G}(U(x,y))=U^{\bot}\widehat{G}(x,y)$. 
When $\mathcal{F}$ defines a figure, with  $U(x,y) \notin F$, we have $\int_{F} G(e^{U(x,y)}) d x d y=0$.
When $d U + d W \in \mathcal{E}^{'(0)}$, this
does not imply  $d U \in \mathcal{E}^{'(0)}$, for instance when the definition of the polar is given 
as linear completion to
$\mathcal{E}^{'(0)}$, corresponding to continuation of the group. Note that $< T(\phi),d U >=<T,dU>(\phi)$,
that is continuation of dU corresponds to the restriction of T, for a fixed $\phi$.

When we only consider sequential movements, $dI \in G$ is reducible. Given d U=d I, 
we have that $d I \notin G_{H}$ ($G_{H}$ of order 1), that is dI does not necessarily have resolution. Given d U=d I implies 
that d V=0,  with  $d V \in G_{H}$, we have presence of resolution.
Assume $d U=pd I$, for p polynomial, given a rectifiable boundary $(1/p) d I \in G$. Assume e d I an 
approximative inverse to p dI, where p is reduced, e can be selected as
symbol to a very regular parametrix, further $e \rightarrow 0$ in $\infty$ (strict condition). The approximation property 
$p \rightarrow 1$ implies that $e \rightarrow 1$ outside traces (modulo terms of negative type), this defines an approximative group.

\newtheorem{prop27}[prop32]{Lemma}
\begin{prop27}
Given d V irreducible =dI, with dV=d WU, for some d W corresponding to a transformation and $d U \in G$, 
we have dW=dI or d U=dI. When $d I \notin G$, we have $dW \notin G$ on the domain. Irreducibles are 
dependent of the order of the group and the topology. Thus, we can simultaneously have $d W=d U^{-1}$ 
in $G_{4}$, given $d I \in G_{4}$.
\end{prop27}

Given $f(w,u)=f_{1}(w)f_{2}(u)$ regular, we have that $f=0$ implies that $f_{1}=0$ or $f_{2}=0$. Given 
regular approximations, we can assume w=0 or u=0. In particular $d I \in G_{2}$ implies that $d I \in G_{4}$.
Assume dV=dWU is irreducible, given contingent regular neighborhoods of singularities relative $G_{2}$, 
they can be approximated by 
$G_{4}$ analogous with a normal tube,
that is $d WU_{j}=\rho_{j} d I$, j=1,2 and $\rho_{j}$ regular.
 
\emph{We assume for a very regular boundary $\{ \Gamma_{i} \}$, existence of w surharmonic, 
with  $\lim_{p \rightarrow q}w(p)=0$ on $\Gamma_{1}$. Given an absolute continuous 
continuation  $\Gamma_{i-1} \rightarrow \Gamma_{1}$,
the boundary has a continuation to a strict condition. We assume for the other segments 
$\Gamma_{i}=\{ u_{i} \}$ $u_{i}$ subharmonic, that $\lim_{p \rightarrow q} \sup_{i} u(p) \leq f(q)$, 
$\forall \{ \Gamma_{i} \}$, where f is real and bounded.}
(cf. \cite{Parreau51},\cite{Collingwood66})

Assume $\Gamma=\{\Gamma_{0},\Gamma_{1} \}$ of bounded variation and $\Omega$ one-sided relative $\Gamma$, 
but not necessarily one-sided relative $\Gamma_{i}$.
Further, $ \phi  : \Gamma_{1} \rightarrow \Gamma_{0}$
absolute continuous. Then we have existence of $\phi^{-1}$. Assume 0 isolated in $\Gamma_{1}$. 
Given a transformation $\Phi$ (we may have $d \Phi \notin G$) 
that regularly traces $\Gamma_{0}$, when $\phi(0)=P \in \Gamma_{0}$, 
with $\phi^{-1}(P)=0$, then $(\phi,\Phi)$ corresponds to a normal tube.

Consider convergence for $\mid \mid V \mid - \mid I \mid \mid \leq \delta$, that corresponds to $\mid I \mid - \delta < \mid V \mid < \mid I \mid + \delta$. 
That is we consider $\mid \int f d \mid V \mid - \int f d \mid I \mid \mid \leq \delta$.
Assume $d \mid V \mid=d V_{+} - d V_{-}$, where we assume $ d \mid V \mid$ of bounded variation. Thus, given $d \mid I \mid=0$, 
we have $d I_{+} \sim d I_{-}$,
that is a two-sided limit. A sufficient condition for regularity for $d \mid V \mid \rightarrow d \mid I \mid$ is thus 
an approximation property and a two-sided limit.

Given f analytic in a bi-disk, for a branch of f, we have $(u,v) = (u,\varphi(u))$,
with  $\varphi$ continuous in $u$.
Consider $G^{\bot}$ as a continuation of G, where f is regular on G for sequential movements, 
but irregular on traces, then
we have that  on $R(G,G^{\bot})$ (considered as a discontinuous extended group), index is well-defined. Given $G^{\bot}$ are linearly independent of G,
a genuine continuation of G is defined. 

\newtheorem{prop29}[prop32]{Lemma}
\begin{prop29}
Assume f=P/Q, then we can define conjugation according to $Q Uf=VP$. Note that $\mid Qf \mid <R$ implies
$\mid f \mid < R' $, for constants R,R', when Q is reduced. Given $0 \neq Q$ has no invariant sets, Q=$V^{-1} U Q$, where
$d V^{-1} U=d I \in G$ implies that Q=0, that is Q(u,v)=0 over a trace set, that defines a polar to f.
\end{prop29}

Note that $If=\int f d I=\int f'(t) d t$, that is convergence in dI over $f \in \mathcal{E}^{0}$,
corresponds to convergence in dt for $f \in \mathcal{E}^{1}$. 
Convergence for I is relative
$\mid f- f_{0} \mid$, convergence for d I is relative $\mid f' - f_{0}' \mid$.

Assume $d U(\phi^{j})/\phi^{j} \rightarrow 0$, when $j > N$, that is invariants only exist for 
$j \leq N$. Given restriction defines the ideal, it is maximal for $j \leq N$, that is invariants can 
be given as a zero-set to an annihilator $W^{\bot}$. However, we do not necessarily have $d W^{\bot} \in G$.
Note that d U=d I trivial does not imply that U=I is trivial. In this case it is for the action of U necessary 
to discuss a polar set.

Assume F=const over $\Delta =\{ u=v \}$, defines multivalentness over $\Delta$. 
A trace approximation between
leaves, can be defined as the limit of $d U=\rho d U_{1},d V=\sigma d U_{2}$ and $\rho/\sigma \rightarrow 1$
through $d U_{1}=d U_{2}$.
Retraction of formulas is then not well-defined.
$\Delta$ defines a carrier (``porteur'') for the limit. $\underline{H} \leq f(q)$ means that when u,f real that f estimates the exponential type 
for $e^{u}$. Given f=0 on bd $\Omega$, $e^{u}$ is locally bounded on bd $\Omega$. Given a very regular boundary, 
we can assume $\underline{H} \leq f \leq \overline{H}$,
further $\underline{H}=\overline{H}$, gives a set characteristic
for solvability. (cf. \cite{Parreau51}).

\subsection{Traces}

Note that f=g in C does not imply that f=g in H. 
Further, $Uf=\int f d I=\int f'(t) d t=f=\int f(u) d u=\int f(0) d u$, that is we assume $\int d u=1$.
Then, Uf=f iff f is absolute continuous, for instance on a domain $\Sigma$, where $f \in L_{ac}^{1}(\Sigma)$. 
Note that $<UT,\phi>=<T,\phi>$ iff $<T,{}^{t}U \phi >=<T,\phi>$, that is invariants in H corresponds
using transponation to invariants in $\mathcal{E}^{' (0)}$.

Consider d U + d V=d I - d W, corresponding to a projective decomposition where d W=0, 
that is a dW is planar (in (u,v)). 
Consider $\Sigma (d U_{j} + d V_{j})=d I -\Sigma d W_{j}$, where $d W_{j}=0$
is dependent of conjugation. The left hand side, has discontinuous convergence and 
$\Sigma + \Sigma d W_{j}$ has continuous convergence.
Given $d W_{j} \geq 0$, we have that $\Sigma d W_{j}=0$ implies that $d W_{j}=0$ $\forall j$, 
that is $\cap_{j} \{ d  W_{j}=0 \} = \{ 0 \}$ generates the space.

\newtheorem{prop13}[prop32]{Lemma}
\begin{prop13}
Assume $f \in L^{1} \cap C^{2}$, with $f(u-v) \equiv 0$ on traces, that are continued to $\infty$. Then we have $f_{u}=-f_{v}$ on traces. 
Further, $f_{vu}=-f_{vv}$. Given a strict condition, we can thus determine index relative $(G_{2})_{H}$ 
for the restriction to $L^{1}_{ac}$ and for f=f(u,v) over the trace set.
\end{prop13}

Consider the range R(dU), with the inner of the range of a subgroup removed, for instance
 $R(d U) \cap R(d V)^{\bot}$, where $R(d V)^{\bot}=R(d V^{\bot})$
(not necessarily unique), that is ``holes'' that are defined by traces. Consider $d U \rightarrow d V \rightarrow d V^{\bot}$
where both arrows are conjugation, but possibly of different nature (not unique). Given 
$d V^{\bot} \rightarrow d I$ regularly,
under preservation of dimension, then the dimension can be determined for R(dV). 

Assume instead the removed sets are defined by $d I \notin (G_{2})_{H}$ (planar). Given a bounded semi-algebraic set $\{ F < b \}$ in $K(G_{4})$,
the boundary can be defined by F=b, thus $F=\Pi  F_{j}$ and the boundary can be decomposed into algebraic sets of
different order. The corresponding sub-level set in $G_{2}$, is unbounded. Example: $F=pf_{0}$, where $p=b$ with $f_{0}$ continuous and locally =1.
Example: $F=P/Q$, 
given Q reduced, according to Nullstellensatz,
$\{ P \leq b \} \subset \{ P/Q \leq b' \}$. When $F \rightarrow 0$ in $\infty$, we have $\infty \in \{ F \leq b \}$.
Given $P \prec \prec Q$, we have $\{ Q \leq b \} \subset \{ P \leq b \}$. 

Assume $\mathcal{F}$ a hyperboloid and consider movements characterized by $U \mathcal{F}=\mathcal{F}$ iff 
$U^{\bot} \mathcal{F}=\mathcal{F}$. 
For euclidean scaling, we have $U \mathcal{F} \cap \mathcal{F}=\{ P \}$ (discrete). Note that when U is non-trivial 
and $P \in \mathcal{F}$, we have $UP \notin \mathcal{F}$.
Assume $\mathcal{F}$ is defined by $ R(d U) \cap R(d V^{\bot})$. Starting from $d U \rightarrow dV \rightarrow  d V^{\bot}$, 
the trace set is defined by $dU=dV^{\bot}$. 
Reflexivity corresponds to $d V=d U^{\bot}$. That is rotation of reflexive traces can be used to generate $\mathcal{F}$. 

When the trace set generates a figure $\mathcal{F}$ in $H_{2}$ (2-dimensional web), given $\mathcal{F}$ a rotational surface, the symmetry is with respect to $d U_{1}$,
that is $f(u_{1},u_{2},u_{3})=f(u_{2},u_{3},u_{1})=const$ on $\mathcal{F}$, in particular $U_{1} \rightarrow (U_{2},U_{3})$ is reversible. 
$U_{1}U_{2}=U_{2}U_{1}$ corresponds to a normal action, outside the trace.

Consider the condition $f(u)=f(w)$, according to $W^{-1}Uf=f$. Given $\Phi$ is a sub-nuclear transformation, 
such that $\Phi \rightarrow I$ preserves subnuclearity, we get a nuclear model, for instance a projective decomposition, 
where $d W \rightarrow 0$ regularly. Sufficient for this, is that $f,\delta f / \delta w \rightarrow 0$ regularly.
Assume $d I \in (G_{2})_{H}(\Omega)$, where $\Omega$ a domain of holomorphy. Given (x,y,z) an inner point, 
such that (u,v) defines a neighborhood in $\Omega$. On bd $\Omega$, f can thus not be continued 
by $H(G_{2})$. However it can be continuable by $C(G_{2})$. 
Given $Uf \in H$, with  $f \in C^{0}$, under a strict condition, we can assume $Uf \in L^{1}$ for $f \in C^{0} \cap L^{1}$.
Thus, we have existence of $f \in C^{0} \cap L^{1}$, that can not be continued regularly outside a domain of holomorphy $\Omega$, 
by $U^{-1}$, that is we assume $d I \notin G_{H}(\Omega)$, when $d U \in G_{L^{1}}(\Omega)$. The problem is analogous 
to the problem of reducing G to translation (and rotation).

\subsection{Reduction of the group}

A decomposition in linear factors, gives a multi-linear kernel. X irreducible has linear convex neighborhoods, 
for instance traces are approximated by sequential movements. X is symplectic, if for instance 
$d U^{\bot}(f)=0$ implies that $d U(f)=d I(f)$ and $d U(f)=0$, that is f=0. A projective decomposition, implies that the trace set is trivial.
Spirals give a non-regular approximation property, in particular it is not a normal model.

\newtheorem{prop3}[prop32]{Lemma}
\begin{prop3}
When J is an irreducible component, we have ${}^{\circ} J^{\circ} \simeq J$,
in particular we have existence of $\omega'$, with  $\omega'(u,v)=0$. 
Thus, there is a decomposition 
in linear factors, that is
$\forall (u,v) \in J$, $< \omega_{1}',(u,v)> \ldots < \omega_{k}',(u,v)>=0$, or 
$< \omega_{1}' \otimes \ldots  \otimes \omega_{k}',(u,v)>=0$,
that defines a maximal ideal. 
\end{prop3}

A resolution of dI ($=d V^{-1}U$) over an irreducible in $G_{2}$, implies that $d I \notin G_{2}$. 
Taking the union over traces in $K(G_{4})$, d I can in the plane, be approximated by $(K(G_{4}))_{H}$.
dU=dV irreducible in $G_{\mathcal{H}}$,
is interpreted as G has a topology compatible with  $\mathcal{H}$, for instance $UT(\phi)=T({}^{t}U \phi)$ and 
$UT(\phi)=VT(\phi)$ for $T \in \mathcal{H}$ and $\forall \phi \in \mathcal{H}'$, implies that ${}^{t}(U-V) \phi=0$, implies that U=V.
Starting from $H(\Omega)$, with topology for Hausdorff uniformities, where $\Omega$ is very regular relative bd $\Omega$, we can construct 
$\Gamma$ very regular and so large that all given singularities can be approximated by $\Omega$.
When J is given by $\omega(u,v)=0$, for a transformation $\omega$, given that $\omega \rightarrow dI$ 
regularly, we have that
$\omega$ is locally 1-1. $d I \in G_{H}(\overline{\Omega} \backslash \Gamma)$, for $\Gamma$ very regular
identifies $\Gamma$, but does not identify individual segments.

Example: consider $\rho(u) \rightarrow \rho(x,y)$ regular, for (x,y) on lines, under preservation of dimension. 
Invariants are taken relative H.
Assume linearly independent movements, d U=d I iff u=0, then $N(u,v,w,z)=N(u) \cup \ldots \cup N(z)$
without common zero's, given existence of regular approximations in for instance u,
the zero space is isolated points on a line.  

We define $\mathcal{H} \supset \mathcal{H}_{\alpha}=\lim_{PC,\beta < \alpha} \mathcal{H}_{\beta}$ 
(PC=point-wise convergence), 
that is we assume presence of algebraic approximations, but not absence of non-regular approximations. 
The condition f=df=0, defines singularities for analyticity,
$\delta f / \delta u=0$, but $\delta f / \delta v \neq 0$, corresponds to sequential approximations 
and to regular approximations of higher order in G. 

Assume $V:C^{\infty} \rightarrow C^{\infty}$, that is point-wise  convergence. When V has a regular
approximation property, as $V \rightarrow I$, dimension is preserved.
When v is monotonously regular in $(x_{1},x_{2})$, v does not lower the dimension.
Note that algebraic continuation preserves hypoellipticity. Starting from $C(G_{2}) \supset H(G_{8})$,
when $f_{j} \in H(G_{2})$ considered on a poly-disk $\mid u_{j} \mid  \leq 1$, we have $\lim_{PC} f_{j} \in H(K(G_{4}))$, 
simultaneously $f_{j} \in C(G_{2})$ implies that $\lim_{PC} f_{j} \in \mathcal{H}_{1}(G_{2})$.

A necessary condition for pseudo convexity, is a boundary of order 0, that is when the boundary 
is given by $d U,d V \in G_{2}$, 
we assume $(u,v) \rightarrow (x,y,z)$ is continuous according to a regular approximation. 
When (u,v) is not surjective on (x,y,z), we consider a boundary of higher order
(for instance $K(G_{4})$) or consider a weaker topology than H.

\newtheorem{prop40}[prop32]{Proposition}
\begin{prop40}
Given (I) is maximal, we have $(I) \subset (J)$, with  $U(J) \subset (J)$, for $d U \in G$, implies that(J)=(I). 
Consider $\phi  : (I) \rightarrow (J)$ sub-nuclear. Maximality
implies that Id is sub-nuclear, that is a nuclear space.
\end{prop40}
(cf. \cite{Schwartz58})
The condition $U (J) \subset (J)$ is necessary for a convolution algebra (measure algebra), that is $d U \in G_{\mathcal{E}^{'(0)}}$. 
Given (I) maximal, there is no interpolation property under U. 
Example: (J)=V(I), then $U(J) \subset (J)$ implies that V=I.
Thus given dU generates G, the continuations must leave G, an interpolation property motivates presence of a transformation,
with $d V \notin G$, that is a discontinuous continuation.

\newtheorem{prop42}[prop32]{Lemma}
\begin{prop42}
Given f can be reduced to algebraic in (u,v) at the boundary $\Gamma$, we have zero's on (u,v) at $\Gamma$. 
Given absence of an approximation property for (u,v), consider a very regular boundary, that is
$\Gamma=\{ \Gamma_{i} \}_{i=1}^{4}$ rectifiable, with  $\Gamma_{i}=\{ (u,v)_{i} \}$ and 
$\Gamma_{i} \rightarrow \Gamma_{i+1}$ absolute continuous. 
We assume a strict condition on some segment. In this case, f can be continued to having zero's of measure zero at $\Gamma$.
\end{prop42}

Given index such that $\mathcal{F}=\{ f=c,df=0 \}$, we have that the character of $\mathcal{F}$ determines the conjugations 
and traces that generate the web.
Consider f(u,v) $\rightarrow$ f(x,y,z), assume f symmetric such that $\lim_{u \rightarrow 0}=\lim_{v \rightarrow 0}$, 
where we assume $u \rightarrow v$ are not linearly dependent 
(a strict condition). Presence of a two-sided limit, is dependent of index. Assume $Wf(u,v)=Wf(v,u)$, 
for instance regularization to a projective decomposition, 
which implies presence of a polar set to (u,v). Example: Wf has negative exponential type,
where f has positive exponential type and $W^{-1}$ does not preserve projectivity.
When the approximation property
is not two-sided, we do not have a convex one-sidedness for regularity relative the boundary.

Consider $\pi  : (u,v) \rightarrow (x,y,z)$, such that $(u=0,\mid v \mid=1)=(x,y,z)$. Assume 
$R(d U,d V) \subset X$, that is $(u,v)(x,y,z)=(x_{1},y_{1},z_{1}) \in X$. Assume presence of a polar, 
with $R(d U,d V)^{\bot} \neq \{ \infty \}$. 
Example: assume $(d U,d V)$ of bounded variation and existence of d W=$\omega d I$, with $R(d W) \subset R(dU,dV)^{\bot}$ 
and $1/\omega=0$ in $\infty$, the extended group satisfies a necessary condition for 
rational convexity for (d U,d V, dW). Example: assume f(u,v)=0 implies that $u \sim 0$, but not $(u,v)=0$, 
for instance when G is such that every d U=d I corresponds to d V=0, with $d V \in G$. 
Assume $d U^{\bot}=\rho d U$, where $\rho$ is regular, when $d U^{\bot}=\sigma d I$,
$d U=\vartheta d I$, we have $\rho = \sigma / \vartheta$, given absence of ess.sng in $\infty$, 
we can select $\sigma,\vartheta$ as holomorphic (\cite{Cousin95}). Pseudo convexity under a strict condition is characterized by an algebraic transversal (cf. \cite{Oka60}). 
Assume $\rho=\sigma / \vartheta \rightarrow 0$, with  $d \rho \neq 0$, we have absence of essential singularities.
Assume for a fixed f, $\{ f=0 \}=\Omega=N(J)$ a domain of holomorphy, then J=J($\rho,f)$ can be selected as minimal, 
($J^{\bot}$ maximal),
note $J \subset rad J \simeq I(\Omega)$. Given $J_{b}=\{ \rho < b \}$, we have $J_{b} \subset J_{b+\epsilon}$,
thus we can determine $b>0$ minimal, such that $\Omega$ is maximal for regularity. 
Example: Assume Wf $\in L^{1}$ (or H), when $f \in C(G_{2})$, for instance given $f d W=\rho f d I$, 
where $\rho \rightarrow 0$ faster than f $\rightarrow \infty$ and $\rho f \in L^{1}$.

\subsection{Resolution in convex components}

Assume a figure is defined using regular transformations $\Phi_{j}$ of translations. The number of traces 
(for instance reflection axes), gives the order for the figure. $\Phi^{\bot}(d U)=d U^{\bot}$ has spiral traces. 
$\overline{\Phi}(d U)=d \overline{U}$ gives ``harmonic'' traces,
given $\Phi^{\diamondsuit}(dU)=d U^{\diamondsuit}$ preserve zero-lines, d U is analytic.  The figure is discussed as the range of
conjugated pairs of movements. The orthogonal to traces includes conjugation relative projective decomposition, 
for which we have two-sided limits. Note that $d U^{\bot}= d V^{\bot}$ does not imply that $d U=d V$.
Consider $R(d U)^{\bot}$, as given by $\{ d V_{j} \}$. Connectivity is dependent of index and of conjugation.
The conjugation is assumed  volume preserving, but not necessarily preserving type.

Assume a regular approximation principle $(u,v) \rightarrow (x,y,z)$ and presence of a separating functional in (u,v), 
then it is necessary for preservation of separating functional, that $u=v$ in a point. When
u=v on a line, the corresponding movement in (x,y,z) does not necessarily have determined tangent, 
thus we have a non-regular approximation property.

Consider one-sidedness relative two axes, such that $R_{1}R_{2} \Omega$ is a convex compact set, where $R_{2}^{-1}$ 
exists. For instance decompose the plane in quadrants, such that $R_{1}R_{2} \Omega \simeq \Omega_{1}^{-} \cap \Omega_{2}^{-} \cup \Omega_{1}^{+} \cap \Omega_{2}^{+}$. 
When the axes can be bijectively identified, this defines a convex compact set in the plane.
Note that when $U \rightarrow U^{\bot}$ is a
contact transform, completion to $L^{1}$ implies possible presence of continuum for $U^{\bot}$.

Assume $R(G \times G)^{\bot}$ a manifold and that L divides this manifold into $\{ \Omega_{j} \} =\{ R(d W_{j} \}$).
A bifurcation can be given by $\Omega_{j}=\Omega_{j}^{+} \cup \Omega_{j}^{-}$, that is $\Omega_{j}$ is not simply 
connected. Every measure of bounded variation $dW=\rho dI$, has a decomposition in components with a determined tangent.
For instance, $d (U \times  V)^{\bot}=d W^{+} - d W^{-}$. Envelop of $d W=dI$, defines a very regular boundary $\Gamma=\{ \Gamma^{\pm} \}$,
given an regular approximation principle for one of the segments.
That is $d \mid W \mid=d \mid I \mid$, does not require
$d W=d I$.

$d I \in G_{H}$ gives a decomposition in convex irreducibles.
For instance consider B: $d U - d I > 0$, with  bd B : d U=d I,
where we note that linearly independent movements have separate invariants.
Example: $dU=dI$ iff dV=0, 
further dV=dI iff dW=0, that is dU+dV=dI is a conjugation relative projective decomposition. 
Assume $d W=\rho d U$, where d U=d I 
is a rectifiable boundary. Consider in particular W(U,V), where d W/d V=0 on bd B.

Assume the complement to R(dU) has a decomposition in linear factors, $R(d U)^{\bot} \simeq \{ d V_{j} \}$. 
In the simplest case, $d U \rightarrow d I$ implies that $d V_{j} \rightarrow 0$ $\forall j$. 
Assume $d V_{1} \rightarrow d V_{2}$ continuous,
but not sub-nuclear, then we have presence of non-projective models.
Example: $d V_{2} \rightarrow d V_{1}$ does not preserve bounded variation. 
Consider $F(\gamma) \rightarrow \gamma \in H$ and $UF(\gamma) \rightarrow 
{}^{t} U \gamma \in H$ (regularization), where we regard $\gamma$ as a pseudo base to a corresponding ideal of holomorphy (\cite{Oka60}). 
Given F harmonic
relative (u,v), there is a pseudo base relative (u,v), 
(for instance a projective decomposition).
Note that $F(\gamma)(u,v) \rightarrow f(x,y,z)$, as $(u,v) \rightarrow (x,y,z)$, assumes an approximation property.

Assume $U \rightarrow I \rightarrow V$ defines analytic continuation, in a neighborhood of u=0. 
Example: $f(u,v)=g(u^{\bot},v)$ regular, where $U^{\bot}=U-I$, that is $u^{\bot}=v$,
 gives a projective decomposition for f. 
$d U^{\bot}=\rho d (U-I)$ preserves zero-lines. The equation defines $U^{\bot}$, given d I=0 has 
a resolution of zero-lines. For $G_{H}$, analytic sets have a decomposition in irreducibles. 
d I $\in G$ gives a decomposition in one-parameter groups.

\newtheorem{prop22}[prop32]{Lemma}
\begin{prop22}
Consider $<d U^{\bot},d U>=0$,
through topological completion to $L^{1}$. Consider $\sigma \sim <d U^{\bot},d V> - <d U,d V^{\bot}>$,
where dV is conjugated to dU.
Assume $\mu=d V / d V^{\bot}$, with  $1/\mu \rightarrow 0$ in $\infty$. Then,
$<\frac{d T}{d V}, <d U^{\bot}, d V>> = <\frac{d T}{d V} d V, d U^{\bot}>=<d T, d U^{\bot}>$.
When $\mu$ is algebraic, we can assume $\phi  \in \mathcal{H}'$ iff $\mu \phi  \in \mathcal{H}'$. 
Further, $<\mu \frac{d T}{d V}, <d U,d V^{\bot}>>=
<d T,d U>$, that is given $d U^{\bot} / d U = \mu = d V / d V^{\bot}$, we have a symplectic form.
For instance $d U^{\bot} d V^{\bot} \simeq  d U d V$.
\end{prop22}

When we use $\sigma \sim 
<\mu d U,d V> - <dU, dV/\mu >$, we have $1/\sigma \rightarrow \pm 0$, as $1/\mu \rightarrow 0$
(and $\rightarrow \infty$), thus $\sigma d I$ has bounded sub-level surfaces $\{ \lambda_{1} < \sigma < \lambda_{2} \}$, 
for scalars $\lambda_{1},\lambda_{2}$, that is we have 
a reduced measure.

Assume $J=\overline{\{ \rho=0 \}} \simeq \{ \rho=const \}$. Note that isolated singularities gives that
J is closed. Thus,
$J=\{ \rho \quad \sigma(\rho)=0 \quad \sigma \in J^{\circ}\}$.
Assume $d V=\rho d U$, with $\rho \in C^{\infty}_{c}$ non-real and such that $\sigma(\rho)=0$, for 
an annihilator $\sigma$, that is $\rho \in (I)$ for a non-real ideal (I).
Assume linearly independent movements, $d V=\rho d U=\rho \nu d I$, with  $\sigma(\rho \nu)=0$
$\forall d V \in G$. In particular, given $\nu \rho_{j}$ quasi-orthogonal polynomials (\cite{Riesz23}), 
such that $\nu \rho_{j}/\mu \in (I)$, for $\mu \in (I)$, we have zeros in a half-space. For instance, $d V_{1} \overline{d V_{2}} - \overline{d V_{1}} d V_{2} \neq 0$,
gives in this case singularities in a half-space. Assume $\mathcal{F}$ is defined by $\sigma(\rho)=\rho$, with  $d U=\rho d I$,
where for instance $\mathcal{F}$ defines the boundary.  
When singularities are characterized by
$d UV=0$ and $d VU \neq 0$, that is $0=d UV= \rho d VU$ implies that $\rho=0$, $\sigma(\rho)=0$ 
implies that $\rho=0$, that is injectivity for $\sigma$.

Assume $\{ u=v \} \rightarrow Hyp$ defines a regular (\cite{Brelot51}) surface. Absence of an approximation property, 
that is $d I \notin (G_{2})_{H}$ but $d I \in (K(G_{4}))_{H}$, is seen as resolution of traces. That is, 
$(u,v,w)$ does now not define traces, but a reducible surface.
For instance, $d V^{-1}U=d I$ implies that $d W=0$, that is $dW \neq d I$ on non-trivial traces, further 
$d (V^{-1}U - I) \rightarrow d W$ defines a projective mapping.

\subsection{Interpolation}
Assume $f=P/Q$ for polynomials P,Q, with $sng f=\{ Q=0 \}$ and sng 1/f = $\{ P=0 \}$, that is $f \rightarrow 1/f$ absolute continuous,
then there is an interpolation analogous to the moment problem (\cite{Riesz23}).

Assume G defines single valued and bounded domains (u,v), where $d U,d V \in G(\Omega)$.
Consider d U/d V=$\rho$ regular, as a condition on convexity with respect to rational functions, 
in particular projective representations according to Radon-Nikodym,
define rationally convex domains. $\rho \rightarrow 0$ somewhere is necessary $(\infty \in \Omega$) for 
integrability.
We thus have domains of holomorphy, that is a boundary of order 0. 

Assume L a line segment and $L \cap \{ F < b \} \neq \emptyset$, that is a convex subset of the 
sub-level surface. Consider $\tilde{L}$ a continuation
to an entire line. Assume $x,y \in \{ F < b \}$, then $x \rightarrow 1/x$ preserves $\tilde{L}$, 
in the same manner for y. Consider $0 \rightarrow x \rightarrow y \rightarrow 0$ 
and $\infty \rightarrow 1/x \rightarrow 1/y \rightarrow \infty$.
Given the duality is linear, lines are preserved.

Example: $f_{1} \leq g \leq f_{2}$, with  $f_{j} \in H(G_{2})$ and $g \in C(G_{2})$.
Example: $g(x,y,1/x,1/y) \rightarrow g(x+1/x,y+1/y)$. When $(x + 1/x,y+1/y)$ is dense in the plane,
we have the same behavior, when $x,y \rightarrow 0$, as when
$x,y \rightarrow \infty$, that is the domain to g is one-sided relative $x \rightarrow 1/x,y \rightarrow 1/y$.

\newtheorem{prop4}[prop32]{Lemma}
\begin{prop4}
Assume  $P_{1} \leq F \leq P_{2}$, where $P_{j}$ are polynomials. Then, there are close to regular 
approximations, non-regular approximations F, $P_{1}(u) \leq F(u,v) \leq P_{2}(u)$. 
The dependence of v for F, means that index for F and $P_{j}$ can be different.
Further, for parametrices according to $P_{j}E_{j}-I=K_{j}$ and $K_{1} \sim K_{2}$, we have that the condition 
$I-E_{2} \leq I-E \leq I-E_{1}$ does not exclude presence of non-trivial kernel for E, 
given $I \neq E_{j}$
\end{prop4}

Assume $d I \in G_{H}$ of bounded variation with separate single -valentness in (u,v), that is a decomposition in linear factors. 
When $d I \notin G$, we have planar multivalentness. Every multivalentness, with holomorphic leaves,
has spiral approximations (\cite{Oka60}). In particular, presence of traces $u_{1}=u_{3}$ implies presence of spiral 
approximating traces.

Consider $L=\{ L_{j} \}$, $j=1,2$, with  $L_{1} \rightarrow L_{2}$ continuous. When L is very regular
and $L_{j}$, j=1,2 are defined by functions of bounded variation, we can select the mapping as absolute continuous. 
Conversely, when $L \subset \{ dW=0 \}$, for a linear transformation dW,
then we have that  Cl(L) includes $L_{1} \rightarrow L_{2}$ continuous. 

Note that dimension assumes $C^{\infty}$ convergence, in particular point-wise. Example: d U $\in G_{\mathcal{H}_{0}}$ 
(\cite{Hartogs37}) and $d W \in G_{\mathcal{H}_{1}}$ implies existence of $d V \in G_{\mathcal{H}_{1}}$, 
with  $d W \leq d V \leq d U$.
Note that spirals do not have dimension, that is point-wise convergence separately in (d U,d V)
does not imply point-wise convergence, spirals can be seen as non-point-wise interpolation of sequential movements.

Example: consider $\lim_{PC} \sup u(p) \leq f(q)$, 
that is $u,f \in \mathcal{H}_{\alpha}$,
thus there is an interpolation property. Note that $F(\gamma) \rightarrow \gamma$ in this case, 
defines a weak interpolation property, that is we have existence of a measure $d V_{F} \in G_{\mathcal{H}_{1}}$
such that $\gamma$ integrable with respect to $d V_{F}$ and $VF(\gamma)$ approximates f(q). 

Consider $F(\gamma) \rightarrow \gamma$ continuous, such that
${}^{t} U \gamma \in H$, then an approximation property gives a deformation to H.
Assume F corresponds to a functional (H', $\mathcal{D}_{L^{1}}'$), for instance ${}^{t} U \gamma \in L^{1} \cap C^{\infty}$, with strict condition. 
On a very regular domain, when f(q) bounded and continuous initial condition to a Dirichlet problem,
we have approximative solutions (\cite{Parreau51}).

Example: (d U,d V) with $d U = \rho d I$ and $d V=(1/\rho) d I$, gives $\rho^{2} d V=d U$, note 
$\rho^{2}$ polynomials, does not imply that $\rho$ is polynomial.
Given $\rho$ is polynomial, there is a decomposition in finitely many linear 
factors.

Example: Assume $\Gamma=\{ \Gamma_{0},\Gamma_{1} \}$ a
very regular boundary and $\lim u=0$ on $\Gamma_{1}$, $u \leq f(q)$ and $\lim v=0$ on $\Gamma_{0}$, with  
$f(q) \leq v$. Given a projective decomposition, $\mid u \mid=1$ on $\Gamma_{0}$.

Assume  dA=0 determines two functionals $d U_{1},d U_{2}$, 
for instance $dA=d U_{1} + d U_{2}$ and $E_{0}=\{ dU_{1}=0 \}$ ,
that is $E_{0} \subset C$ and  $d A \varphi=dI \varphi$, for some element outside $E_{0}$ (cf. \cite{Riesz56}).
The condition for $E_{0}$ dense in C, is given by $M_{A} \geq 1/d$, 
where d denotes the distance between $E_{0}$ and $\varphi$ such that $dA \varphi=dI \varphi$ 
and $M_{A}$ denotes the upper limit for action of dA over $E_{0}$. 
Using (\cite{Lie91}, kap 19,20) we have that a projective decomposition can be given starting from a
boundary condition in (x,y).

Consider $\beta_{1} d U + \beta_{2} d U^{\bot}=0$, with  $\beta_{1} : \beta_{2}$ non-constant in (x,y,z).
For instance $\beta_{1}=d N / d U$ and $\beta_{2}=d N / d U^{\bot}$, we then have d N=0 and 
the condition means
$d U^{\bot} /d U$ non-constant. Then N can be given using an integration in x,y.

\subsection{Subgroups}
 
Starting from uniformities on X (\cite{Bourbaki89}), when $\Delta$ is the diagonal in $X \times X$,
$\cap nbhd \Delta=\Delta$ iff X Hausdorff relative uniformities. Note that for every nbhd $\Delta$=
$\Omega$, we have existence of W, such that
$(x,y) \in \Omega$, with  $(x,z) \in W$ and $(z,y) \in W$, for some $z \in \Omega$. 
Given $\Delta^{c}$ irreducible in $\Omega \backslash \Delta$ and $W \subsetneq \Omega \backslash \Delta$, 
we have that $W=0$. 

Assume $(\Gamma_{0},\Gamma_{1})$ a very regular boundary, with  $\Gamma_{1} \rightarrow \Gamma_{0}$ absolute continuous.
Given Hausdorff topology, that is sng is given by $\Gamma_{0}$, we can consider $\Gamma_{1} \sim \Gamma_{0}$ 
(a deformation), such that arbitrarily close to $\Gamma_{0}$ on $\Gamma_{1} \rightarrow \Gamma_{0}$, there are regular points.
Note that $\Gamma_{0}(G)=\{ \Gamma_{0}^{j} \}$, where not all segments are mapped onto the web.

Given irreducibles in $H_{1}$ (that is traces that result in a 1-dimensional web),
we have removable traces in the plane, that is u=v is mapped on
a line, removable in the plane as constant surfaces in H. Simultaneously, irreducibles u=v can be mapped on $H_{1}$, 
considered as irregular in 3-space. Consider $\Gamma_{1}=\{ \mid u-v \mid =R \}$ and 
$\Gamma_{0}=\{ u=v \}$, with  $\Gamma_{1} \rightarrow \Gamma_{0}$ absolute continuous, then traces can 
in $\Gamma$ be seen as in $H_{0}$,
for a very regular boundary (\cite{Parreau51}). Homotopic deformation is assumed to preserve dimension. 
In the same manner, the homeomorphism $(\Gamma_{0},u_{2}) \rightarrow \mathcal{F}$ is assumed to preserve dimension 
(profile curve). 
  
Note that $\Omega(T(f))=\{ T(f) < \lambda \}$ semi-algebraic does not imply that $\Omega(f)$ is 
semi-algebraic, $\Omega(T(f)) \subset \Omega(f)$, that is $\Omega(f)$ includes a semi-algebraic set.
Assume $\{ P/Q \leq \lambda \}$ with  P,Q polynomials, given Q is reduced, using Nullstellensatz (\cite{Oka60})
$\{ P \leq \lambda' \} \subset \{ P/Q \leq \lambda \}$, for scalars $\lambda,\lambda'$. Assume existence of a polynomial R, 
with  $R \leq P/Q$. 
Then through the condition
of slow growth $QR \leq P$, for Q reduced, we have that $R \leq P$.

Assume singularities are given by f=df=0. 
Example: $F=e^{\phi}$, given $d \phi$ is bounded, where F=0 , we have dF=0 and $d \phi  \sim 0/0$.
$d F(u-v) \equiv 0$, gives $\frac{\delta F}{\delta u} \equiv -\frac{\delta F}{\delta v}$, thus $\frac{\delta^{2} F}{\delta u \delta v}=-\frac{\delta^{2 F}}{\delta v^{2}}$ ($\neq 0$), 
that is we have index. Note that when A is generated by a very regular group, $\int_{A} f'(t) d t=\int_{A} f d I=\int_{bd A} f$.
Assume U preserves singularities, we then have Ud F=UF=0 over sng, that is $< d F,d U >=< F,d U >=0$. 
Example: With these conditions, $\Gamma \subset N(d U)$ implies that $UF=0$ on $\Gamma$,
further f(u)=d f(u)=0.
 
Consider $(x,y,z) \rightarrow (u,v) \rightarrow (n,n^{\bot})$, as a boundary operator, where the last
mapping is continuously reversible. Assume $(n,n^{\bot})$
irreducible, then $d N(U,U^{\bot})$ allows traces on a disk neighborhood. That is, we do not have 
a strict condition and we have a non-normal extension.
The condition $f(n,n^{\bot})=0$ implies that $f(n)=0$ or $f(n^{\bot})=0$, can be interpreted as 
f being irreducible for $(n,n^{\bot})$, that is
given f is irreducible, we have that
$f \mid (n,n^{\bot})$ implies that $f \mid n$ or $f \mid n^{\bot}$. Assume $(d U,0) \in \mathcal{E}^{'(0)}$ 
and continuation through
$d U^{\bot}$ has a representation in $\mathcal{E}'$, with  branch-point analogous to Hartogs. 
Example: $d N(f)=0$ is analytic with $d (N,N^{\bot})(f) \neq 0$
and index indicates presence of traces. 
Simultaneously, $f(n+n^{\bot})=0$ for $n^{\bot} \neq 0$ gives an approximation property.

Consider $F(\gamma) \rightarrow \gamma$ continuous, where $\gamma$ are for instance pseudo bases. 
Consider further
$\{ F < b \} \rightarrow U \gamma$, depending on an inverse lifting principle. That is, consider $f(u,v)=F(\gamma)(u,v)$,
and $\Delta_{u,v}f=0$, that gives a domain for analytic inversion. 
Starting from $UF(\gamma)=F({}^{t} U \gamma)$,
given an approximation property and assuming UF nuclear, the topology can be selected as nuclear. 
Given $UF \in H$, given an approximation property over $\gamma$, we can determine ${}^{t} U^{-1} \gamma$. 
Starting from $d U \in (G_{2})_{C}$,, 
there is d V $(\in K(G_{4}))_{H}$ close to d U, that preserves analyticity
and we have an inverse lifting principle, over dV. 

Consider $C^{\infty} \subset C^{0}$, where dimension is not assumed to be well-defined in $C^{0}$, 
we assume that movements preserve regularity.
Example: Assume $Wf \in C^{\infty}$, for $f \in C^{0}$, with $W^{-1}W=I$ in $C^{\infty}$, 
that is $d W^{-1} \notin G_{\mathcal{E}^{(0)'}}$.
Example: $G_{\mathcal{E}^{(0)'}} \ni  d U \rightarrow d U^{\bot} \in G_{\mathcal{E}'}$, 
where the inverse can be given in $\mathcal{E}'$,
that is completion gives solvability in $\mathcal{E}'$, corresponding to a discontinuous group. 

Consider $Uf \in H$, when $f \in C^{0}$, that is $d U^{-1} \notin G_{H}$. 
Analytic functionals H', have representation in $\mathcal{E}^{'(0)}$ (\cite{Martineau}).
Assume $\rho d I \in \mathcal{E}^{' (0)}$, for instance $\rho \in C_{0}^{\infty}$, consider $ C_{0}^{\infty} \rightarrow C^{\infty}$,
for instance over a pseudo convex domain (relatively compact sub-level surfaces), 
we have an approximation property through truncation, when $\rho d I \in \mathcal{E}^{'(0)}$, as $\rho \rightarrow 1$.

Consider $\mathcal{E}' \rightarrow \mathcal{E}^{' (0)}$, for instance $d U=\rho d I$, with
$\rho dI \in \mathcal{E}^{'(0)}$, over $f \in C^{\infty}$. When we discuss approximative solutions E, 
Ef - I of negative type, considered on a symmetric domain, it is
sufficient to consider E of type 0 (\cite{Dahn13}), in particular it is sufficient to consider measures in phase space.
Note that $e^{L()}$ algebraic does not imply that $e^{-L()}$ algebraic, where L is a scalar product. 
Consider a majorization condition, for instance $\mid \delta f / \delta v \mid \leq \mid f(w) \mid$, 
that is existence of w, such that $f(w) \rightarrow 0$ implies that $\delta f / \delta v \rightarrow 0$. 
In particular, when G is discontinuous, by extending G with W as above, we have a continuous group.

Example: $(U_{1} + U_{2})f=f(u_{1},u_{2})=f(u_{1},0) + f(0,u_{2})$. 
Given $f \in B_{pp}(u_{1},u_{2})$, we have that
$f = f_{1}(u_{1}) \times f_{2}(u_{2})$. 

 Given $d \mu \in G_{\mathcal{E}^{'(0)}}$ defines a transformation, then 
$d I \in G_{\mathcal{E}^{(0)'}}$ is necessary for invertibility,.
Consider the completion to a Hilbert space $\mathcal{H}$,  
then we do not necessarily have sub-nuclear convergence $d \mu \rightarrow d I$, 
that is $d I \notin G_{\mathcal{H}}$, in particular we can have presence of traces.

Example: Assume $L(f)=\int f d \mu$. Consider a transformation $\tilde{L}(f)=\int \tilde{f} d \nu$, 
where $d \mu \rightarrow d \nu$ and further $d \mu_{2} \rightarrow d \nu^{2}$ is Legendre. 
Necessary for $d \nu \in G_{\mathcal{E}^{'(0)}}$ is a strict condition. Assume further for the 
conjugated $d \nu$, that $e^{I(d \nu)} \in G_{\mathcal{E}^{'(0)}}$, with $I^{2}(d \nu)=I(d \nu^{2})$ 
and corresponding to $U^{2}f=\int f d \mu_{2}$. In this manner convergence for $e^{I(d \mu)}$ can be
derived in a dual sense.

Assume $d u \rightarrow d \tilde{u}$ continuation to closedness and $d \tilde{u} \rightarrow d v$ 
conjugation, further $d \tilde{u}(d v)=d u (d \tilde{v})$. The condition $\int d ud v=0$, 
using the volume preserving property
remains constant for continuation along lines. 
Assume $d T(u)=\mid d u \mid = d \mid u \mid$ with T absolute continuous,
then we have $\mid \int d u d v \mid \leq \int dT(u) d T(v) =\int d T(u,v)$.
Assume bdA is given by $T(u,v) =1$ and A is given by $ T(u,v) < 1$. Consider inner points that 
are defined by u=0 and $\mid v \mid < 1$ or $\mid u \mid < 1$ and $\mid v \mid=1$, 
that is a boundary to a bi-disk.
Assume this boundary is domain for regularity for f, then in particular f is regular on A. 

Example: $d U=\rho d U_{1}$,$d V=\sigma d U_{2}$, 
that is $<f, d U \times d V>=<\sigma \rho f, d U_{1} \times d U_{2} >$.
When $d U(f)/f=1$, we have for d U($f^{j})$, that d U($e^{f}$) is not convergent. 

Assume $d U=\rho d V$, then dV=0 implies that $d U=0$. Given $\rho(v)$ is algebraic, there is a zero
close to 0. Consider dU as a continuation of dV. Given a projective decomposition, where f is absolute 
continuous, Uf + Vf=f, when outside the domain for absolute continuous convergence,
$d U + d V \notin G_{H}$, we consider $\in G_{C}$. Thus, given $d U + d V=d I \notin G_{H}$, 
we must use a discontinuous group $ G_{C}$.

\newtheorem{prop25}[prop32]{Lemma}
\begin{prop25}
Assume $F(\gamma)=UF(V \gamma)$, then UF can be determined in $\mathcal{E}'$, given that V 
preserves analyticity and $\gamma \in H$. Assume $\Gamma_{1}$ with isolated singularities, 
$\Gamma_{0}$ with continuum for $f \in H$.
In particular, UF can be defined in $\mathcal{E}^{' (0)}$ over $\Gamma_{1}$, the corresponding 
representation is not in $H'$  over $\Gamma_{0}$. Assume
$\psi  : \Gamma_{1} \rightarrow \Gamma_{0}$ absolute continuous, such that $\psi(f) \in C^{0}$, 
then a measure $\in \mathcal{E}^{'(0)}$ can be
continued to $(\Gamma_{0},\Gamma_{1})$ very regular.
\end{prop25}

\cite{Garding87}
\bibliographystyle{amsplain}
\bibliography{index_2024}

\end{document}